# Directional Differentiability of the Metric Projection in Bochner Spaces

Jinlu Li

Department of Mathematics
Shawnee State University
Portsmouth, Ohio 45662
USA

## Abstract

In this paper, we consider the directional differentiability of metric projection and its properties in uniformly convex and uniformly smooth Bochner space $L_p(S; X)$, in which $(S, \mathcal{A}, \mu)$ is a positive measure space and $X$ is a uniformly convex and uniformly smooth Banach space. Let (arbitrary) $A \in \mathcal{A}$ with $\mu(A) > 0$ and define a subspace $L_p(A; X)$ of $L_p(S; X)$, which is considered as a closed and convex subset of $L_p(S; X)$. We first study the properties of the normalized duality mapping in $L_p(S; X)$ and in $L_p(A; X)$. For any $c \in L_p(A; X)$ and $r > 0$, we define a closed ball $B_A(c; r)$ in $L_p(A; X)$ and a cylinder $C_A(c; r)$ in $L_p(S; X)$ with base $B_A(c; r)$. Then, we investigate some optimal properties of the corresponding metric projections $P_{L_p(A;X)}$, $P_{B_A(c;r)}$ and $P_{C_A(c;r)}$ that include the inverse images, the directional differentiability and the precise solutions of their directional derivatives.

**Mathematics subject classification** (**2010**): 65K10, 90C25, 90C26, 90C48

## 1. Introduction

In this paper, we always, otherwise stated, consider the (standard) metric projection in uniformly convex and uniformly smooth Banach spaces (see [1, 12, 16, 23]). Let $(Z, \lVert \cdot \rVert)$ be a real uniformly convex and uniformly smooth Banach space with topological dual space $(Z^*, \lVert \cdot \rVert_*)$. Let $C$ be a nonempty closed and convex subset of $Z$. Let $P_C: Z \to C$ denote the (standard) metric projection operator, which is a well-defined single-valued mapping, such that, for any $x \in Z$, we have $P_C x \in C$ satisfying
$$\lVert x - P_C x \rVert \le \lVert x - z \rVert, \text{ for all } z \in C.$$

$P_C x$ is called the metric projection of point $x$ onto $C$. $P_C x$ is considered as the best approximation of $x$ by elements of $C$, which is the closest point from $x$ to $C$. For any $y \in C$, the inverse image of $y$ by the metric projection $P_C$ in $Z$ is defined by
$$P_C^{-1}(y) = \{x \in X: P_C(x) = y\}.$$

In particular, if $Z$ is a Hilbert space, then $P_C$ has the following properties:

(i)     The basic variational principle: for any $x \in Z$ and $u \in C$,

$$u = P_C x \quad \Longleftrightarrow \quad \langle x - u, u - z \rangle \ge 0, \text{ for all } z \in C.$$

(ii)     $P_C$ is nonexpansive:

$$\|P_C x - P_C y\| \le \|x - y\|, \text{ for any } x, y \in Z.$$

With help of the above properties, the directional differentiability of $P_C$ in Hilbert spaces have been studied by many authors (see [11, 13, 14,18, 22, 24]).

However, the metric projection $P_C$ in a uniformly convex and uniformly smooth Banach space $Z$ does not enjoy the above simple basic variational principle (i) and the nonexpansive property (ii), in general. Let $J_Z$: $Z \to Z^*$ be the normalized duality mapping in this uniformly convex and uniformly smooth Banach space $Z$, which is a well-defined single-valued mapping (see [1, 29]). With the help of $J_Z$ in uniformly convex and uniformly smooth Banach space $Z$, $P_C$ has the basic variational principle: for any $x \in Z$ and $u \in C$,

$$u = P_C x \quad \Longleftrightarrow \quad \langle J_Z(x - u), u - z \rangle \ge 0, \text{ for all } z \in C.$$

Since $J_Z$ is not a linear operator, it substantially increases the difficulty in studying the directional differentiability of $P_C$ in uniformly convex and uniformly smooth Banach spaces. In [2, 6, 17, 25−26, 30], some types of directional differentiability are studied for some projection operators in Banach spaces, which has been applied to approximation theory, convex program problems, optimal control problems, and so forth (see [3, 4, 19, 20]).

Since Bochner spaces can be considered as special cases of Banach spaces, in this paper, we concentrate to study the directional differentiability of the metric projection in uniformly convex and uniformly smooth Bochner spaces. The definition of the directional differentiability follows that in [17]. We have more colorful properties of the directional differentiability of the metric projection in Bochner spaces, which are extension of the results studied in [17] in uniformly convex and uniformly smooth Banach spaces.

This paper is organized as follows. At first, we review some concepts and properties of uniformly convex and uniformly smooth Bochner spaces. Then, we study the properties and the representations of the normalized duality mapping in some subspaces of Bochner spaces, which will be used to study the properties of the metric projection in Bochner spaces (see [7, 15]). Let $C$ be a closed and convex subset of a uniformly convex and uniformly smooth Bochner space. In section 4, we will investigate some properties and solutions of $P_C$. In section 5, following the definition in [17], we give the definition of the directional differentiability of the metric projection operator and investigate some properties. Then we find the precise representations of the directional derivatives of $P_C$. In sections 4 and 5, we especially consider some special cases of $C$: proper subspaces, closed balls and closed cylinders of uniformly convex and uniformly smooth Bochner spaces. As applications of the results obtained in sections 4 and 5, in section 6, we study the properties and solutions of the directional derivatives of $P_C$ in Hilbert spaces.

## 2. Preliminaries

### 2.1. Uniformly convex and uniformly smooth Bochner spaces

In this section, we review some concepts and properties of uniformly convex and uniformly smooth Bochner spaces (see [5, 7, 9, 10, 15, 21, 27, 28]). Let $(S, \mathcal{A}, \mu)$ be a positive measure space. Let $(X, \|\cdot\|_X)$ be a real uniformly convex and uniformly smooth Banach space with topological dual space $(X^*, \|\cdot\|_{X^*})$. Let $\langle \cdot, \cdot \rangle$ denote the real canonical evaluation pairing between $X^*$ and $X$. For any $A \in \mathcal{A}$ and $x \in X$, let $1_A \otimes x$ denote the $X$-valued simple function on $S$ defined by

$$(1_A \otimes x)(\text{s}) = 1_A(s) \otimes x = \begin{cases} x, & \text{if } s \in A, \\ 0, & \text{if } s \notin A, \end{cases} \qquad \text{for any } s \in S,$$

where $1_A$ denotes the characteristic function of $A$ on the space $S$. For an arbitrary given positive integer $n$, let $\{A_1, A_2, \ldots A_n\}$ be a finite collection of mutually disjoint subsets in $\mathcal{A}$ with $0 < \mu(A_i) < \infty$, for all $i = 1, 2, \ldots n$. Let $\{x_1, x_2, \ldots x_n\} \subseteq X$ and let $\{a_1, a_2, \ldots a_n\}$ be a set of real numbers. Then, $\sum_{i=1}^{n} a_i(1_{A_i} \otimes x_i)$ is called a $\mu$-simple function from $S$ to $X$ (See Definition 1.1.13 in [15]).

Throughout this paper, we take positive numbers $p$ and $q$, with $1 < p, q < \infty$ satisfying $\frac{1}{p} + \frac{1}{q} = 1$. Let $(L_p(S; X), \|\cdot\|_{L_p(S;X)})$ be the Lebesgue-Bochner function space that is called the Bochner space based on the measure space $S$ and the Banach space $X$, which is a real uniformly convex and uniformly smooth Banach space. More precisely speaking, $L_p(S; X)$ is the Banach space of $\mu$-equivalent classes of strongly measurable functions $f\colon S \to X$ with norm:

$$\|f\|_{L_p(S;X)} = \left(\int_S \|f(s)\|_X^p d\mu(s)\right)^{\frac{1}{p}} < \infty, \text{ for } 1 < p < \infty.$$

The dual space of $(L_p(S; X), \|\cdot\|_{L_p(S;X)})$ is $(L_q(S; X^*), \|\cdot\|_{L_q(S;X^*)})$. In this paper, we let $\langle \cdot, \ \cdot \rangle_{L_p}$ denote the real canonical evaluation pairing between the uniformly convex and uniformly smooth Bochner spaces $L_q(S; X^*)$ and $L_p(S; X)$. For easy referee, we list some properties of Bochner integrals and Bochner spaces below.

(B$_1$). $\left\|\int_S f d\mu\right\|_X \leq \int_S \|f\|_X d\mu$, for every $f \in L_p(S, X)$;

(B$_2$). Every $\varphi \in L_q(S, X^*)$ defines a bounded linear functional $\varphi \in (L_p(S; X))^*$ by the formula

$$\langle \varphi, f \rangle_{L_p} = \int_S \langle \varphi(\omega), f(\omega) \rangle d\mu(\omega), \text{ for every } f \in L_p(S, X).$$

It satisfies

$$\|\varphi\|_{(L_p(S; X))^*} = \|\varphi\|_{L_q(S; X^*)}, \text{ for any } \varphi \in L_q(S, X^*).$$

(B$_3$). $L_2(S; X)$ is a Hilbert space $\iff$ $X$ is a Hilbert space.

(B$_4$). For an arbitrary $A \in \mathcal{A}$ with $0 < \mu(A) < \infty$ and for any $x, y \in X$, we have

(i) $\frac{1}{\mu(A)^{\frac{1}{p}}}(1_A \otimes x) \in L_p(S; X)$;

(ii) $\left\|\frac{1}{\mu(A)^{\frac{1}{p}}}(1_A \otimes x) \pm \frac{1}{\mu(A)^{\frac{1}{p}}}(1_A \otimes y)\right\|_{L_p(S;X)} = \|x \pm y\|_X$;

(iii) The mapping $x \to \frac{1}{\mu(A)^{\frac{1}{p}}}(1_A \otimes x)$ (isometrically) embeds $X$ into $L_p(S; X)$.

See [5, 9, 15, 21, 27] for more properties and more concepts of Bochner spaces.

For the considered uniformly convex and uniformly smooth Banach space $X$, the normalized duality mapping $J_X\colon X \to X^*$ is a single-valued mapping satisfying

$$\langle J_X x, x \rangle = \|J_X x\|_{X^*}\|x\|_X = \|x\|_X^2 = \|J_X x\|_{X^*}^2, \text{ for any } x \in X.$$

The normalized duality mapping $J_X$ in uniformly convex and uniformly smooth Banach space $X$ has many useful properties (see [1, 29] for more details). For example,

(i) $J_X: X \to X^*$ is one to one and onto;

(ii) $J_X$ is a continuous;

(iii) $J_X$ is a homogeneous operator;

(iv) $J_X$ is uniformly continuous on each bounded subset of $X$.

In particular, for the real Banach space $l_p$ and $L_p(S)$, with $1 < p < \infty$, the normalized duality mapping holds the following analytic representations.

(a) For any $x = (t_1, t_1, \dots) \in l_p$ with $x \neq \theta$,

$$(J_{l_p} x)_n = \frac{|x_n|^{p-1} \text{sign}(x_n)}{\|x\|_{l_p}^{p-2}} = \frac{|x_n|^{p-2} x_n}{\|x\|_{l_p}^{p-2}}, \text{ for } n = 1, 2, \dots.$$

(b) For any $h \in L_p(S)$ with $h \neq \theta$,

$$(J_{L_p(S)} h)(s) = \frac{\|h(s)\|_X^{p-1} \text{sign}(h(s))}{\|h\|_{L_p(S)}^{p-2}} = \frac{\|h(s)\|_X^{p-2} h(s)}{\|h\|_{L_p(S)}^{p-2}}, \text{ for all } s \in S.$$

The normalized duality mappings in the uniformly convex and uniformly smooth Bochner space $L_p(S; X)$ is denoted by $J_{L_p(S,X)}$, which is abbreviated as $J_p$ if there is no confusion caused. Then, by the properties of normalized duality mappings, both $J_X$ and $J_p$ are single valued, one to one and onto continuous maps. The normalized duality mappings in the dual spaces $X^*$ and $L_q(S; X^*)$ are respectively denoted by $J_{X^*}$ and $J_q^*$, if there is no confusion caused. They have the following properties and analytic representations, which are proved in [7].

**Corollary 3.2 in [7]**. *Let $A \in \mathcal{A}$ with $0 < \mu(A) < \infty$. Then, for any $x \in X$ with $x \neq 0$, we have*

$$J_p(1_A \otimes x)(s) = \frac{\mu(A)^{\frac{1}{p}}}{\mu(A)^{\frac{1}{q}}} (1_A \otimes J_X x)(s), \text{ for all } s \in S.$$

*It is equivalent to*

$$J_p(\frac{1}{\mu(A)^{\frac{1}{p}}} (1_A \otimes x))(s) = \frac{1}{\mu(A)^{\frac{1}{q}}} (1_A \otimes J_X x)(s), \text{ for all } s \in S.$$

**Corollary 3.3 in [7]**. *$J_p$ maps every $\mu$-simple function in $L_p(S; X)$ to $\mu$-simple function in $L_q(S; X^*)$ with respect to the same partition in $S$. Moreover, for an arbitrary given $\mu$-simple function $\sum_{i=1}^n (1_{A_i} \otimes x_i)$ in $L_p(S; X)$, we have*

$$J_p\big(\sum_{i=1}^n (1_{A_i} \otimes x_i)\big)(s) = \frac{1}{\big(\sum_{j=1}^n \|x_j\|_X^p \mu(A_j)\big)^{\frac{1}{q} - \frac{1}{p}}} \sum_{i=1}^n \|x_i\|_X^{p-2} (1_{A_i} \otimes J_X x_i)(s), \text{ for all } s \in S.$$

**Corollary 3.4 in [7]**. *For any $f \in L_p(S; X)$, let $\{f_n\}$ be a sequence of $\mu$-simple functions in $L_p(S; X)$ satisfying*

$$f_n \to f, \text{ in } L_p(S; X), \text{ as } n \to \infty.$$

*Then $\{J_p f_n\}$ is a sequence of $\mu$-simple functions in $L_q(S; X^*)$ such that*

$$J_p f_n \to J_p f, \text{ in } L_q(S; X^*), \text{ as } n \to \infty.$$

## 2.2 The function of smoothness of uniformly convex and uniformly smooth of Banach spaces

For a uniformly convex and uniformly smooth Banach space $(X, \|\cdot\|)$ with topological dual space $(X^*, \|\cdot\|_*)$, let $S(X)$ be the unit sphere of $X$, that is, $S(X) = \{x \in X : \|x\| = 1\}$. Then, it is well-known that $X$ is uniformly smooth if and only if the limit

$$\lim_{t \downarrow 0} \frac{\|x+tv\| - \|x\|}{t},$$

exists uniformly for all $(x, v) \in S(X) \times S(X)$. Then, in [17], we introduced the following definition.

**Definition 2.1 in [17].** Let $X$ be a uniformly convex and uniformly smooth Banach space. Define $\psi_X : S(X) \times S(X) \to \mathbb{R}$ by

$$\psi_X(x, v) = \lim_{t \downarrow 0} \frac{\|x+tv\| - \|x\|}{t}, \text{ for any } (x, v) \in S(X) \times S(X). \tag{2.1}$$

$\psi_X$ is called the function of smoothness of this Banach space $X$. Since $X$ is uniformly convex and uniformly smooth, then, the limit (2.1) is attained to $\psi_X(x, v)$ uniformly for $(x, v) \in S(X) \times S(X)$ (see Takahashi [26]).

For the convenience and simplicity, we extend the concept of function of smoothness of uniformly smooth Banach space $X$ from $S(X) \times S(X)$ to $X \times (X \setminus \{\theta\})$.

**Definition 2.1.** Let $X$ be a uniformly convex and uniformly smooth Banach space. Define $\Psi_X : X \times (X \setminus \{\theta\}) \to \mathbb{R}$ by

$$\Psi_X(x, v) = \lim_{t \downarrow 0} \frac{\|x+tv\| - \|x\|}{t}, \text{ for any } (x, v) \in X \times (X \setminus \{\theta\}). \tag{2.2}$$

$\Psi_X$ is called the (extended) function of smoothness of this Banach space $X$. Next, we show that $\Psi_X$ is indeed an extension of $\psi_X$ defined in (2.1) from $S(X) \times S(X)$ to $X \times (X \setminus \{\theta\})$.

**Lemma 2.2.** Let $X$ be a uniformly convex and uniformly smooth Banach space. Then, for any $(x, v) \in X \times (X \setminus \{\theta\})$, we have

$$\Psi_X(x, v) = \begin{cases} \|v\|, & \text{for } x = \theta, \\ \|v\| \psi_X \left( \frac{x}{\|x\|}, \frac{v}{\|v\|} \right), & \text{for } x \neq \theta. \end{cases} \tag{2.3}$$

In particular,

$$\Psi_X(x, x) = \|x\|, \text{ for any } x \neq \theta.$$

*Proof.* If $x = \theta$, then it is clear to have

$$\Psi_X(\theta, v) = \|v\|, \text{ for any } v \neq \theta.$$

For any $(x, v) \in (X \setminus \{\theta\}) \times (X \setminus \{\theta\})$, we have

$$\Psi_X(x, v) = \lim_{t \downarrow 0} \frac{\|x+tv\| - \|x\|}{t}$$
$$= \lim_{t \downarrow 0} \frac{\|x\| \left( \left\| \frac{x}{\|x\|} + \frac{\|v\|t}{\|x\|} \frac{v}{\|v\|} \right\| - \left\| \frac{x}{\|x\|} \right\| \right)}{t}$$

$$= \lim_{t \downarrow 0} \frac{\|v\|}{\|x\|} \frac{\|x\| \left( \left\| \frac{x}{\|x\|} + \frac{\|v\| t}{\|x\|} \frac{v}{\|v\|} \right\| - \left\| \frac{x}{\|x\|} \right\| \right)}{\frac{\|v\| t}{\|x\|}}$$

$$= \|v\| \psi_X \left( \frac{x}{\|x\|}, \frac{v}{\|v\|} \right). \qquad \square$$

Since in the last section of this paper, we will study the properties of the metric projection in Hilbert spaces, so, in next lemma, we consider the function of smoothness of Hilbert spaces, which are considered as special cases of uniformly convex and uniformly smooth Banach spaces.

**Lemma 2.3.** *Let* $(H, \|\cdot\|)$ *be a Hilbert space with inner product* $\langle \cdot, \ \cdot \rangle$. *Then, for any* $(x, v) \in H \times (H \backslash \{\theta\})$, *we have*

$$\Psi_H(x, v) = \begin{cases} \|v\|, & \text{if } x = \theta, \\ \frac{1}{\|x\|} \langle x, v \rangle, & \text{if } x \neq \theta. \end{cases} \qquad (2.4)$$

*Proof.* The proof of this lemma is straight forward and it is omitted here. $\qquad \square$

**Notation**s: Let $L_p(S; X)$ be a uniformly convex and uniformly smooth Bochner space with $1 < p < \infty$. For the simplicity of notations, the function of smoothness of this uniformly smooth Bochner space $L_p(S; X)$ is rewritten as:

$$\Psi_{L_p(S,X)} \equiv \Psi_p.$$

## 3. The normalized duality mapping in subspaces of Bochner spaces

**Definition 3.1**. Let $L_p(S; X)$ be a uniformly convex and uniformly smooth Bochner space with $1 < p < \infty$. For an arbitrary $A \in \mathcal{A}$ with $\mu(A) > 0$ and for $f \in L_p(S; X)$, if $f$ satisfies

$$f(s) = \theta, \text{ for } \mu\text{-almost every } s \in S \backslash A,$$

then, $f$ is said to be supported in $A$. We denote the collection of all functionals in $L_p(S; X)$ supported in $A$ as follows

$$L_p(A; X) = \{f \in L_p(S; X): f(s) = \theta, \text{ for } \mu\text{-almost every } s \in S \backslash A\}.$$

We similarly define

$$L_q(A; X^*) = \{\varphi \in L_q(S; X^*): \varphi(s) = \theta, \text{ for } \mu\text{-almost every } s \in S \backslash A\}.$$

**Lamma 3.2**. *$L_p(A; X)$ has the following properties*:

(a) *$L_p(A; X)$ is a subspace of $L_p(S; X)$*;

(b) *If $\mu(S \backslash A) > 0$, then $L_p(A; X)$ is a proper subspace of $L_p(S; X)$*;

(c) *$L_p(A; X)$ is a uniformly convex and uniformly smooth Bochner space with the measure space* $(A, \mathcal{A}|_A, \mu|_A)$ *such that*

$$(L_p(A; X))^* = L_q(A; X^*).$$

*Proof.* The proof of this lemma is straight forward and it is omitted here. $\qquad \square$

Let $A \in \mathcal{A}$ with $\mu(A) > 0$. For $f \in L_p(S; X)$, we write

$$f_A(s) = \begin{cases} f(s), \text{for } s \in A, \\ \theta, \ \text{ for } s \notin A. \end{cases}$$

**Lamma 3.3**. *Let $L_p(S; X)$ be a uniformly convex and uniformly smooth Bochner space with $1 < p < \infty$. Let $A \in \mathcal{A}$ with $\mu(A) > 0$. For any $f \in L_p(S; X)$, we have*

(a) $f_A \in L_p(A; X)$;

(b) $(\lambda f)_A = \lambda f_A$, *for any $\lambda \in \mathbb{R}$;*

(c) $(f + g)_A = f_A + g_A$, *for any $f, g \in L_p(S; X)$.*

*Proof.* The proof of this lemma is straight forward and it is omitted here. □

**Lamma 3.4**. *Let $L_p(S; X)$ be a uniformly convex and uniformly smooth Bochner space with $1 < p < \infty$. For $A \in \mathcal{A}$ with $\mu(A) > 0$ and for $f \in L_p(S; X)$, we have*

(a) $J_p f_A \in L_q(A; X^*)$ *and*

$$(J_p f_A)_A = J_p f_A;$$

(b) *Furthermore,*

$$f \in L_p(A; X) \Longrightarrow J_p f \in L_q(A; X^*).$$

*That is, $J_p \colon L_p(A; X) \longrightarrow L_q(A; X^*)$;*

(c) *Moreover, for any $f \in L_p(S; X)$ with $f_A \neq \theta$, we have*

$$(J_p f)_A = \frac{\|f_A\|_{L_p(A; X)}^{p-2}}{\|f\|_{L_p(S; X)}^{p-2}} J_p f_A.$$

*More precisely, for $s \in A$.*

$$(J_p f)(s) = \frac{\|f_A\|_{L_p(S; X)}^{p-2}}{\|f\|_{L_p(S; X)}^{p-2}} (J_p f_A)(s).$$

(d) *In general, if $\mu(S \backslash A) > 0$, then*

$$(J_p f)_A \neq J_p f_A.$$

*Proof.* Proof of (a). We calculate

$$\langle J_p f_A, f_A \rangle_{L_p}$$
$$= \int_S \langle (J_p f_A)(s), \ f_A(s) \rangle d\mu(s)$$
$$= \int_A \langle (J_p f_A)(s), \ f_A(s) \rangle d\mu(s) + \int_{S \backslash A} \langle (J_p f_A)(s), \ f_A(s) \rangle d\mu(s)$$
$$= \int_A \langle (J_p f_A)(s), \ f_A(s) \rangle d\mu(s)$$
$$= \int_A \langle (J_p f_A)_A(s), \ f_A(s) \rangle d\mu(s) + \int_{S \backslash A} \langle (J_p f_A)_A(s), \ f_A(s) \rangle d\mu(s)$$
$$= \int_S \langle (J_p f_A)_A(s), \ f_A(s) \rangle d\mu(s)$$

$$= \langle (J_p f_A)_A, \; f_A \rangle_{L_p}.$$

By the definition of $J_p$, this implies

$$
\begin{aligned}
&\|J_p f_A\|_{L_q(S;X^*)} \|f_A\|_{L_p(S;X)} \\
&= \langle J_p f_A, \; f_A \rangle_{L_p} \\
&= \langle (J_p f_A)_A, \; f_A \rangle_{L_p} \\
&\leq \|(J_p f_A)_A\|_{L_q(S;X^*)} \|f_A\|_{L_p(S;X)}.
\end{aligned}
\tag{3.1}
$$

It follows that

$$\|J_p f_A\|_{L_q(S;X^*)} \leq \|(J_p f_A)_A\|_{L_q(S;X^*)}.$$

It is clear to see

$$\|(J_p f_A)_A\|_{L_q(S;X^*)} \; \leq \; \|J_p f_A\|_{L_q(S;X^*)}.$$

This implies

$$\|(J_p f_A)_A\|_{L_q(S;X^*)} \; = \; \|J_p f_A\|_{L_q(S;X^*)}.
\tag{3.2}$$

By (3.1), we have

$$\langle (J_p f_A)_A, f_A \rangle_{L_p} = \|(J_p f_A)_A\|_{L_q(S;X^*)} \|f_A\|_{L_p(S;X)}.
\tag{3.3}$$

On the other hand, by (3.2), we have

$$\|f_A\|_{L_p(S;X)} = \|J_p f_A\|_{L_q(S;X^*)} = \|(J_p f_A)_A\|_{L_q(S;X^*)}.
\tag{3.4}$$

(3.3) and (3.4) imply

$$(J_p f_A)_A = J_p f_A.$$

Hence, $J_p f_A \in L_p(A; X).$

It is clear that part (b) follows from part (a) immediately. Moreover, we give a directly proof for part (b) by using the representations of $J_p$. By Proposition 3.1 in [7], for any $f \in L_p(S; X)$ with $f_A \neq \theta$, we have

$$
\begin{aligned}
(J_p f_A)(s) &= \frac{\|f_A(s)\|_X^{p-2} J_X(f_A(s))}{\|f_A\|_{L_p(S;X)}^{p-2}}, \quad \text{for all } s \in S \\
&= \begin{cases} \dfrac{\|f(s)\|_X^{p-2} J_X(f(s))}{\|f_A\|_{L_p(A;X)}^{p-2}}, & \text{for all } s \; \in \; A, \\[2mm] \theta, & \text{for all } s \; \in \; S \backslash A. \end{cases}
\end{aligned}
\tag{3.5}
$$

This implies $J_p f_A \in L_q(A; X^*)$. Next, we prove (c). For any $f \in L_p(S; X)$ with $f_A \neq \theta$, by Proposition 3.1 in [7], we have

$$(J_p f)(s) = \frac{\|f(s)\|_X^{p-2} J_X(f(s))}{\|f\|_{L_p(S;X)}^{p-2}}, \quad \text{for all } s \in S.$$

In particular, for all $s \in A$, we have

$$(J_p f)(s) = \frac{\|f(s)\|_X^{p-2} J_X(f(s))}{\|f\|_{L_p(S;X)}^{p-2}},$$

$$= \frac{\|f_A\|_{L_p(S;X)}^{p-2}}{\|f\|_{L_p(S;X)}^{p-2}} \frac{\|f_A(s)\|_X^{p-2} J_X(f_A(s))}{\|f_A\|_{L_p(S;X)}^{p-2}}$$

$$= \frac{\|f_A\|_{L_p(S;X)}^{p-2}}{\|f\|_{L_p(S;X)}^{p-2}} (J_p f_A)(s), \text{ for all } s \in A. \tag{3.6}$$

This proves (c). In case, if $\mu(S \backslash A) > 0$, then, we can choose $f \in L_p(S; X)$ such that $\|f\|_{L_p(S;X)} \neq \|f_A\|_{L_p(A;X)}$. For such a $f \in L_p(S; X)$, by (3.5) and (3.6), we have

$$(J_p f)(s) \neq (J_p f_A)(s), \text{ for all } s \in A, \text{ at which } f(s) \neq \theta.$$

From this, part (d) follows immediately. □

**Definition 3.5.** Let $(S, \mathcal{A}, \mu)$ be a measure space. Let $\mathbb{N}$ be the set of all positive integers and $M$ a nonempty subset of $\mathbb{N}$. Let $\{A_n \in \mathcal{A} : n \in M\}$ be a family of subsets of $S$ satisfying $\mu(A_n) > 0$, for every $n \in M$. If the following conditions are satisfied

$$A_n \cap A_m = \emptyset, \text{ for } n \neq m \text{ and } \bigcup_{n \in M} A_n = S,$$

then, $\{A_n \in \mathcal{A} : n \in M\}$ is called a strong partition of $S$.

**Proposition 3.6.** *Let $L_p(S; X)$ be a uniformly convex and uniformly smooth Bochner space with $1 < p < \infty$. Let $\{A_n \in \mathcal{A} : n \in M\}$ be a strong partition of $S$. For any $f \in L_p(S; X)$ with $f \neq \theta$, then*

$$J_p f = \frac{1}{\|f\|_{L_p(S;X)}^{p-2}} \sum_{n \in M} \|f_{A_n}\|_{L_p(A_n;X)}^{p-2} J_p f_{A_n}.$$

*Proof.* Let $f \in L_p(S; X)$ with $f \neq \theta$. Define

$$M_f = \{n \in M : \|f_{A_n}\|_{L_p(A_n;X)} \neq 0\}.$$

Since $\{A_n \in \mathcal{A} : n \in M\}$ is a strong partition of $S$, by the representation of $J_p$ given in Proposition 3.1 in [7], for all $s \in S$, we have

$$(J_p f)(s) = \frac{\|f(s)\|_X^{p-2} J_X(f(s))}{\|f\|_{L_p(S;X)}^{p-2}}$$

$$= \sum_{n \in M} \frac{\|f_{A_n}(s)\|_X^{p-2} J_X(f_{A_n}(s))}{\|f\|_{L_p(S;X)}^{p-2}}$$

$$= \sum_{n \in M_f} \frac{\|f_{A_n}(s)\|_X^{p-2} J_X(f_{A_n}(s))}{\|f\|_{L_p(S;X)}^{p-2}}$$

$$= \frac{1}{\|f\|_{L_p(S;X)}^{p-2}} \sum_{n \in M_f} \frac{\|f_{A_n}\|_{L_p(S;X)}^{p-2} \|f(s)\|_X^{p-2} J_X(f_{A_n}(s))}{\|f_{A_n}\|_{L_p(S;X)}^{p-2}}$$

$$= \frac{1}{\|f\|_{L_p(S;X)}^{p-2}} \sum_{n \in M_f} \|f_{A_n}\|_{L_p(S;X)}^{p-2} \frac{\|f_{A_n}(s)\|_X^{p-2} J_X(f_{A_n}(s))}{\|f_{A_n}\|_{L_p(S;X)}^{p-2}}$$

$$= \frac{1}{\|f\|_{L_p(S;X)}^{p-2}} \sum_{n \in M_f} \|f_{A_n}\|_{L_p(S;X)}^{p-2} (J_p f_{A_n})(s)$$

$$= \frac{1}{\|f\|_{L_p(S;X)}^{p-2}} \sum_{n \in M} \|f_{A_n}\|_{L_p(S;X)}^{p-2} (J_p f_{A_n})(s), \text{ for all } s \in S.$$

This implies

$$J_p f = \frac{1}{\|f\|_{L_p(S;X)}^{p-2}} \sum_{n \in M} \|f_{A_n}\|_{L_p(A_n;X)}^{p-2} J_p f_{A_n}. \qquad \square$$

**Corollary 3.7**. *Let $L_p(S; X)$ be a uniformly convex and uniformly smooth Bochner space with $1 < p < \infty$. For any $A \in \mathcal{A}$ with $\mu(A) > 0$ and $\mu(S \backslash A) > 0$, let $L_p(A; X)$ be the proper subspace of $L_p(S; X)$ given in Definition 3.1. For any $f \in L_p(S; X)$ with $f \neq \theta$, we have*

$$J_p f = \frac{1}{\|f\|_{L_p(S;X)}^{p-2}} \left( \|f_A\|_{L_p(A;X)}^{p-2} J_p f_A + \|f_{S \backslash A}\|_{L_p(S \backslash A;X)}^{p-2} J_p f_{S \backslash A} \right).$$

## 4. The metric projection in uniformly convex and uniformly smooth Bochner spaces

### 4.1. The metric projection onto closed subspaces of uniformly convex and uniformly smooth Bochner spaces

Recall that, for $A \in \mathcal{A}$ with $\mu(A) > 0$, $(L_p(A; X), \|\cdot\|_{L_p(A;X)})$ is a subspace of $L_p(S; X)$ given in Definition 3.1. The topological dual space of $L_p(A; X)$ is $L_q(A; X^*)$. If the given subset $A$ of $S$ satisfies $\mu(S \backslash A) > 0$, then $L_p(A; X)$ is a proper subspace of $L_p(S; X)$. In this subsection, we study the metric projection from $L_p(S; X)$ onto $L_p(A; X)$. It helps us to study the directional differentiability of the metric projection in next section.

**Lemma 4.1**. *Let $A \in \mathcal{A}$ with $\mu(A) > 0$ and $\mu(S \backslash A) > 0$. For any $f \in L_p(S; X)$, we have*

(a) $P_{L_p(A;X)}(f) = f_A$;

(b) $P_{L_p(A;X)}(\lambda f) = \lambda f_A$, *for any $\lambda \in \mathbb{R}$;*

*Proof*. Proof of (a). For any $g \in L_p(A; X)$, we calculate

$$\langle J_p(f - f_A), f_A - g \rangle_{L_p}$$
$$= \int_S \langle J_p(f - f_A)(s), f_A(s) - g(s) \rangle d\mu(s)$$
$$= \int_A \langle J_p(f - f_A)(s), f_A(s) - g(s) \rangle d\mu(s) + \int_{S \backslash A} \langle J_p(f - f_A)(s), f_A(s) - g(s) \rangle d\mu(s)$$
$$= \int_A \langle J_p(\theta_A), f_A(s) - g(s) \rangle d\mu(s) + \int_{S \backslash A} \langle J_p(f)(s), \theta \rangle d\mu(s)$$
$$= \int_A \langle \theta_A, f_A(s) - g(s) \rangle d\mu(s) + \int_{S \backslash A} \langle J_p(f)(s), \theta \rangle d\mu(s)$$
$$= 0, \text{ for every } g \in L_p(A; X).$$

By the basic variational principle of $P_{L_p(A;X)}$, this proves (a) of this lemma. Then (b) of this lemma follows from part (a) of this lemma and part (b) of Lemma 3.3. $\square$

**Corollary 4.2**. *Let $A \in \mathcal{A}$ with $\mu(A) > 0$ and $\mu(S \backslash A) > 0$. Then, we have,*

(i) $P_{L_p(A;X)}(f) = \theta_A, \text{ for any } f \in L_p(S \backslash A; X).$

(ii) $P_{L_p(A;X)}^{-1}(\theta_A) = L_p(S \backslash A; X).$

(iii)     *For any $h \in L_p(A; X)$, $P_{L_p(A;X)}^{-1}(h)$ is a closed and convex cone in $L_p(S; X)$ with vertex at $h$ such that*

$$P_{L_p(A;X)}^{-1}(h) = h + L_p(S \backslash A; X).$$

*Proof.* Parts (i) and (ii) of this corollary follow from Lemma 4.1 immediately. We only show part (iii). For any $h \in L_p(A; X)$, by part (a) in Lemma 4.1, we have

$$\begin{aligned}
P_{L_p(A;X)}^{-1}(h) &= \{f \in L_p(S; X) : f_A = h\} \\
&= h + \{f \in L_p(S; X) : f_A = \theta_A\} \\
&= h + P_{L_p(A;X)}^{-1}(\theta_A) \\
&= h + L_p(S \backslash A; X).
\end{aligned}$$ □

It is well-known that, in general, the metric projection operator in uniformly convex and uniformly smooth Banach spaces is not nonexpansive. However, in particular, in uniformly convex and uniformly smooth Bochner spaces, when the considered closed and convex subset $C$ is a proper subspace defined in Definition 3.1, the metric projection operator is nonexpansive, which can be proved by using Lemma 4.1 as follows.

**Corollary 4.3**. *Let $A \in \mathcal{A}$ with $\mu(A) > 0$ and $\mu(S \backslash A) > 0$. Then, $P_{L_p(A;X)}$ is nonexpansive.*

*Proof.* For any $f, g \in L_P(S; X)$, by Lemma 4.1, we have

$$P_{L_p(A;X)}(f) = f_A \text{ and } P_{L_p(A;X)}(g) = g_A.$$

This implies

$$\begin{aligned}
&\left\| P_{L_p(A;X)}(f) - P_{L_p(A;X)}(g) \right\|_{L_p(S;X)} \\
&= \left\| P_{L_p(A;X)}(f) - P_{L_p(A;X)}(g) \right\|_{L_p(A;X)} \\
&= \|f_A - g_A\|_{L_p(A;X)} \\
&\leq \|f - g\|_{L_p(S;X)}.
\end{aligned}$$ □

## 4.2. **The metric projection onto closed balls in subspaces of uniformly convex and uniformly smooth Bochner spaces**

For any $r > 0$ and $v \in L_p(A; X)$, we define the closed, open balls (open in $L_p(A; X)$) and the sphere in $L_p(A; X)$ with radius $r$ and with center at $v$, respectively, by

$$B_A(v, r) = \left\{ f \in L_p(A; X) : \left( \int_A \|f(s) - v(s)\|_X^p d\mu(s) \right)^{\frac{1}{p}} \leq r \right\},$$

$$B_A^o(v, r) = \left\{ f \in L_p(A; X) : \left( \int_A \|f(s) - v(s)\|_X^p d\mu(s) \right)^{\frac{1}{p}} < r \right\},$$

$$S_A(v, r) = \left\{ f \in L_p(A; X) : \left( \int_A \|f(s) - v(s)\|_X^p d\mu(s) \right)^{\frac{1}{p}} = r \right\}.$$

We note that $B_A(v, r)$ is a nonempty closed, bounded and convex subset in $L_p(S; X)$. If $\mu(S \backslash A) > 0$, then $B_A(v, r)$ is not a (closed) ball in $L_p(S; X)$ and $B_A^o(v, r)$ is not open in $L_p(S; X)$. However, $B_A^o(v, r)$ and $S_A(v, r)$ form a partition of $B_A(v, r)$. In particular, if $v = \theta$, then $B_A(\theta, r)$, $B_A^o(\theta, r)$ and $S_A(\theta, r)$ are denoted by $B_A(r)$, $B_A^o(r)$ and $S_A(r)$, respectively. In this subsection, we consider the properties of the metric projection onto closed balls $B_A(r)$. All results about $B_A(r)$ can be analogously extended to $B_A(v, r)$. For $r > 0$, based on $B_A(r)$, we define the following subsets in $L_p(S; X)$:

$$C_A(r) = \{f \in L_p(S; X): f_A \in B_A(r)\},$$

$$C_A^o(r) = \{f \in L_p(S; X): f_A \in B_A^o(r)\}.$$

$C_A(r)$ and $C_A^o(r)$ are called cylinders in $L_p(S; X)$ with bases $B_A(r)$ and $B_A^o(r)$ in $L_p(A; X)$, respectively. If $\mu(S\backslash A) > 0$, then $C_A(r)$ is a closed, unbounded and convex subset in $L_p(S; X)$ and $C_A^o(r)$ is an open, unbounded and convex subset in $L_p(S; X)$. In particular, if $\mu(S\backslash A) = 0$, then $L_p(A; X)$ coincides with $L_p(S; X)$ and, in this case, $B_A(v, r)$, $B_A^o(v, r)$, $S_A(v, r)$, $C_A(r)$ and $C_A^o(r)$ are respectively denoted by $B(v, r)$, $B^o(v, r)$, $S(v, r)$, $C(r)$ and $C^o(r)$. Moreover, in the case that $\mu(S\backslash A) = 0$, we have

$$C(r) = B(r) \quad \text{and} \quad C^o(r) = B^o(r).$$

Next, we calculate the values of the metric projection from $L_p(S; X)$ onto $B_A(r)$.

**Theorem 4.4**. *For $r > 0$ and $g \in L_p(S; X)$, we have*

(a) $P_{B_A(r)}(g) = g$, *for $g \in B_A(r)$;*

(b) $P_{B_A(r)}(g) = \dfrac{r}{\|g\|_{L_p(A; X)}} g$, *for $g \in L_p(A; X)\backslash B_A(r)$;*

(c) $P_{B_A(r)}(g) = g_A$, *for $g \in C_A(r)\backslash L_p(A; X)$;*

(d) $P_{B_A(r)}(g) = \dfrac{r}{\|g_A\|_{L_p(A; X)}} g_A$, *for $g \in L_p(S; X)\backslash(C_A(r) \cup L_p(A; X))$.*

*Proof.* Part (a) is clear. We prove part (b). For any given $g \in L_p(A; X)\backslash B_A(r)$, $g$ must satisfies $\|g\|_{L_p(S; X)} = \|g\|_{L_p(A; X)} > r$. Then, for any $f \in B_A(r)$, we have

$$\langle J_p(g - \frac{r}{\|g\|_{L_p(S; X)}} g), \ \frac{r}{\|g\|_{L_p(S; X)}} g - f\rangle$$

$$= \left(1 - \frac{r}{\|g\|_{L_p(A; X)}}\right) \langle J_p(g), \ \frac{r}{\|g\|_{L_p(A; X)}} g - f\rangle$$

$$= \left(1 - \frac{r}{\|g\|_{L_p(A; X)}}\right) \left(r\|g\|_{L_p(A; X)} - \langle J_p(g), f\rangle\right)$$

$$\geq \left(1 - \frac{r}{\|g\|_{L_p(A; X)}}\right) \left(r\|g\|_{L_p(A; X)} - \|J_p(g)\|_{L_q(A; X^*)} \|f\|_{L_p(A; X)}\right)$$

$$= \left(1 - \frac{r}{\|g\|_{L_p(A; X)}}\right) \left(r\|g\|_{L_p(A; X)} - \|g\|_{L_p(A; X)} \|f\|_{L_p(A; X)}\right)$$

$$= \left(\|g\|_{L_p(A; X)} - r\right) \left(r - \|f\|_{L_p(A; X)}\right)$$

$$\geq 0, \text{ for all } f \in B_A(r).$$

Since $\dfrac{r}{\|g\|_{L_p(S; X)}} g = \dfrac{r}{\|g\|_{L_p(A; X)}} g \in S_A(r)$, by the basic variational principle of $P_{B_A(r)}$, this implies

$$P_{B_A(r)}(g) = \frac{r}{\|g\|_{L_p(A; X)}} g, \quad \text{for } g \in L_p(A; X)\backslash B_A(r).$$

Proof of (c). Since $g \in C_A(r)\backslash L_p(A; X)$, then $g_A \in B_A(r)$ and $g \notin L_p(A; X)$. It follows that

$$\|g_A\|_{L_p(A;\, X)} \le r \quad \text{and} \quad \|g_A\|_{L_p(A;\, X)} < \|g\|_{L_p(S;\, X)}.$$

Then, for any $f \in B_A(r) \subseteq L_p(A; X)$, we have

$$\langle J_p(g - g_A),\ g_A - f \rangle$$

$$= \int_S \langle J_p(g - g_A)(s),\ g_A(s) - f(s) \rangle d\mu(s)$$

$$= \int_A \langle J_p(g - g_A)(s),\ g_A(s) - f(s) \rangle d\mu(s) + \int_{S \setminus A} \langle J_p(g - g_A)(s),\ g_A(s) - f(s) \rangle d\mu(s)$$

$$= \int_A \langle J_p(\theta_A),\ g_A(s) - f(s) \rangle d\mu(s) + \int_{S \setminus A} \langle J_p(g - g_A)(s),\ \theta \rangle d\mu(s)$$

$$= \int_A \langle \theta_A,\ g_A(s) - f(s) \rangle d\mu(s) + \int_{S \setminus A} \langle J_p(g - g_A)(s),\ \theta \rangle d\mu(s)$$

$$= 0, \text{ for every } f \in B_A(r) \subseteq L_p(A; X).$$

Since $g_A \in B_A(r)$, by the basic variational principle of $P_{B_A(r)}$, this implies

$$P_{B_A(r)}(g) = g_A, \quad \text{for } g \in C_A(r) \backslash L_p(A; X).$$

The proof of part (d) is similar to the proof of part (b). For $g \in L_p(S; X) \backslash (C_A(r) \cup L_p(A; X))$, $g \notin C_A(r) \cup L_p(A; X)$, then $g_A \notin B_A(r)$. By $g \notin L_p(A; X)$, it follows that

$$\|g_A\|_{L_p(A;\, X)} > r \quad \text{and} \quad \|g_A\|_{L_p(A;\, X)} < \|g\|_{L_p(S;\, X)}.$$

Then, $f \in B_A(r)$ implies $\frac{r}{\|g_A\|_{L_p(S;\, X)}} g_A - f \in L_p(A; X)$. By Lemma 3.4 and Corollary 3.7, we have

$$\left\langle J_p\left(g - \frac{r}{\|g_A\|_{L_p(S;\, X)}} g_A\right),\ \frac{r}{\|g_A\|_{L_p(S;\, X)}} g_A - f \right\rangle$$

$$= \left\langle J_p\left(g_{S \setminus A} + \left(1 - \frac{r}{\|g_A\|_{L_p(A;\, X)}}\right) g_A\right),\ \frac{r}{\|g_A\|_{L_p(S;\, X)}} g_A - f \right\rangle$$

$$= \left(1 - \frac{r}{\|g_A\|_{L_p(A;\, X)}}\right) \left\langle J_p(g_A),\ \frac{r}{\|g_A\|_{L_p(A;\, X)}} g_A - f \right\rangle$$

$$= \left(1 - \frac{r}{\|g_A\|_{L_p(A;\, X)}}\right) \left(r\|g_A\|_{L_p(A;\, X)} - \langle J_p(g_A), f \rangle\right)$$

$$\ge \left(\|g_A\|_{L_p(A;\, X)} - r\right)\left(r - \|f\|_{L_p(A;\, X)}\right)$$

$$\ge 0, \text{ for all } f \in B_A(r).$$

By $\|g_A\|_{L_p(A;\, X)} > r$, we have $\frac{r}{\|g_A\|_{L_p(S;\, X)}} g_A = \frac{r}{\|g_A\|_{L_p(A;\, X)}} g_A \in B_A(r)$, by the basic variational principle of $P_{B_A(r)}$, this implies

$$P_{B_A(r)}(g) = \frac{r}{\|g_A\|_{L_p(S;\, X)}} g_A, \quad \text{for } g \in L_p(S; X) \backslash (C_A(r) \cup L_p(A; X)). \qquad \square$$

Notice that if $\mu(S \backslash A) = 0$, then $L_p(A; X)$ coincides with $L_p(S; X)$. In this case, we have

$$B_A(v, r) = B(v, r) \quad \text{and} \quad C_A(r) = B(r).$$

Then, by Theorem 4.4, we have

**Corollary 4.5**. *For $r > 0$ and $g \in L_p(S; X)$, we have*

(a) $P_{B(r)}(g) = g$, for $g \in B(r)$;

(b) $P_{B(r)}(g) = \dfrac{r}{\|g\|_{L_p(S;\,X)}}\, g$, for $g \in L_p(S; X) \backslash B(r)$.

By applying Theorem 4.4 to simple functions in $L_p(S; X)$, we have the following results.

**Corollary 4.6**. *Let $A \in \mathcal{A}$ with $\mu(A) = 1$ and $x \in X$. For $r > 0$, we have*

(a) $P_{B_A(r)}(1_A \otimes x) = 1_A \otimes x$, if $\|x\|_X \leq r$;

(b) $P_{B_A(r)}(1_A \otimes x) = \dfrac{r}{\|x\|_X}(1_A \otimes x)$, if $\|x\|_X > r$.

*Proof.* By $\mu(A) = 1$, we calculate

$$\|1_A \otimes x\|_{L_p(S;\,X)} = \|x\|_X.$$

Then, this corollary follows from Corollary 4.5 immediately. □

### 4.3. The metric projection onto closed cylinders in uniformly convex and uniformly smooth Bochner spaces

In next theorem, we calculate the values of the metric projection from $L_p(S; X)$ onto a closed and convex cylinder $C_A(r)$ in $L_p(S; X)$.

**Theorem 4.7**. *For any $r > 0$ and $g \in L_p(S; X)$, we have*

(a) $P_{C_A(r)}(g) = g$, for $g \in C_A(r)$;

(b) $P_{C_A(r)}(g) = \dfrac{r}{\|g_A\|_{L_p(A;\,X)}}\, g_A + g_{S \backslash A}$, for $g \in L_p(S; X) \backslash C_A(r)$.

*Proof.* Part (a) is clear. We prove part (b). For any given $g \in L_p(S; X) \backslash C_A(r)$, it satisfies

$$\|g_A\|_{L_p(S;\,X)} = \|g_A\|_{L_p(A;\,X)} > r. \tag{4.1}$$

By part (c) in Lamma 3.3, we notice that

$$\left( \frac{r}{\|g_A\|_{L_p(A;\,X)}}\, g_A + g_{S \backslash A} \right)_A = \frac{r}{\|g_A\|_{L_p(A;\,X)}}\, g_A.$$

By (4.1), this implies

$$\frac{r}{\|g_A\|_{L_p(A;\,X)}}\, g_A + g_{S \backslash A} \in C_p(r). \tag{4.2}$$

It is clear that, for $g \in L_p(S; X)$, we have the following decomposition

$$g = g_A + g_{S \backslash A}. \tag{4.3}$$

Then, for any $f \in C_A(r)$, by $J_p(g_A) \in L_q(A; X^*)$ and (4.1), we have

$$\langle J_p \left( g - \left( \frac{r}{\|g_A\|_{L_p(A;\,X)}} g_A + g_{S\setminus A} \right) \right), \left( \frac{r}{\|g_A\|_{L_p(A;\,X)}} g_A + g_{S\setminus A} \right) - f \rangle$$

$$= \langle J_p \left( (g_A + g_{S\setminus A}) - \frac{r}{\|g_A\|_{L_p(A;\,X)}} g_A - g_{S\setminus A} \right), \left( \frac{r}{\|g_A\|_{L_p(A;\,X)}} g_A + g_{S\setminus A} \right) - (f_A + f_{S\setminus A}) \rangle$$

$$= \left( 1 - \frac{r}{\|g_A\|_{L_p(A;\,X)}} \right) \langle J_p(g_A), \frac{r}{\|g_A\|_{L_p(A;\,X)}} g_A - f_A \rangle + \left( 1 - \frac{r}{\|g_A\|_{L_p(A;\,X)}} \right) \langle J_p(g_A), \ g_{S\setminus A} - f_{S\setminus A} \rangle$$

$$= \left( 1 - \frac{r}{\|g_A\|_{L_p(A;\,X)}} \right) \langle J_p(g_A), \frac{r}{\|g_A\|_{L_p(A;\,X)}} g_A - f_A \rangle$$

$$= \left( 1 - \frac{r}{\|g_A\|_{L_p(A;\,X)}} \right) \frac{r}{\|g_A\|_{L_p(A;\,X)}} \|g_A\|_{L_p(A;\,X)}^2 - \left( 1 - \frac{r}{\|g_A\|_{L_p(A;\,X)}} \right) \langle J_p(g_A), f_A \rangle$$

$$\geq \left( 1 - \frac{r}{\|g_A\|_{L_p(A;\,X)}} \right) \left( r \|g_A\|_{L_p(A;\,X)} - \|g_A\|_{L_p(A;\,X)} \|f_A\|_{L_p(A;\,X)} \right)$$

$$= \left( \|g_A\|_{L_p(A;\,X)} - r \right) \left( r - \|f_A\|_{L_p(A;\,X)} \right)$$

$$\geq 0, \text{ for all } f \in C_A(r).$$

By (4.2) and by the basic variational principle of $P_{C_A(r)}$, this implies

$$P_{C_A(r)}(g) = \frac{r}{\|g_A\|_{L_p(A;\,X)}} g_A + g_{S\setminus A}, \quad \text{for } g \in L_p(S; X) \setminus C_A(r). \qquad \square$$

**Remarks on Corollary 4.5 and Theorem 4.7.** As we noticed before Corollary 4.5, if $\mu(S\setminus A) = 0$, then $L_p(A; X)$ coincides with $L_p(S; X)$. In this case, we have $C_A(r) = B(r)$ and

$$g_{S\setminus A} = \theta, \text{ for } g \in L_p(S; X).$$

Then, we can see that Corollary 4.5 can be induced by Theorem 4.7.

## 5. The directional differentiability of the metric projection in uniformly convex and uniformly smooth Bochner spaces

### 5.1. The directional differentiability of the metric projection in uniformly convex and uniformly smooth Banach spaces

In [17], the Gâteaux directional differentiability of $P_C$ in uniformly convex and uniformly smooth Banach spaces is defined (see Definition 4.1 in [17]). This topic was deeply studied in Hilbert spaces in [18]. In this section, we recall the definitions and properties obtained in [17].

In this subsection, let $(Z, \|\cdot\|)$ be a real uniformly convex and uniformly smooth Banach space with topological dual space $(Z^*, \|\cdot\|_*)$. Let $C$ be a nonempty closed and convex subset of $Z$. For an arbitrary $x \in Z$ and for a vector $v \in Z$ with $v \neq \theta$, the directional derivative of $P_C$ at point $x$ along direction $v$ is defined by

$$P_C'(x; v) = \lim_{t \downarrow 0} \frac{P_C(x+tv) - P_C(x)}{t},$$

provided the existence of this limit. Then, $P_C$ is said to be (Gâteaux) directionally differentiable and $P_C'(x; v)$ is called the (Gâteaux) directional derivative of $P_C$ at point $x$ along direction $v$. Vector $v$ is called a (Gâteaux) differentiable direction of $P_C$ at $x$. If $P_C$ is (Gâteaux) directionally differentiable at point $x \in Z$ along every direction $v \in Z$ with $v \neq \theta$, then $P_C$ is said to be (Gâteaux) directionally differentiable at point $x \in Z$, which is denoted by

$$P'_C(x)(v) = \lim_{t \downarrow 0} \frac{P_C(x+tv) - P_C(x)}{t}, \text{ for } v \in Z \text{ with } v \neq \theta. \tag{5.1}$$

$P'_C(x)(v)$ is called the (Gâteaux) directional derivative of $P_C$ at point $x$ along direction $v$. Let $D$ be a subset in $Z$. If $P_C$ is (Gâteaux) directionally differentiable at every point $x \in D$, then $P_C$ is said to be (Gâteaux) directionally differentiable on $D \subseteq Z$. Many properties of the (Gâteaux) directional derivative of $P_C$ are proved in [17]. We list some of them below for easy referee.

1. *The following statements are equivalent*

(i)     *$P_C$ is directionally differentiable on $Z$ such that, for every point $x \in Z$,*

$$P'_C(x)(v) = \theta, \text{ for any } v \in Z \text{ with } v \neq \theta;$$

(ii)     *$P_C$ is a constant operator; that is, $C$ is a singleton.*

2. *For every point $x \in Z$, there is least one differentiable direction of $P_C$ at $x$.*

3. *For $x \in Z$ and for $v \in Z$ with $v \neq \theta$, if $P_C$ is directionally differentiable at $x$ along direction $v$, then, for any $\lambda > 0$, $P_C$ is directionally differentiable at $x$ along direction $\lambda v$ such that*

$$P'_C(x; \lambda v) = \lambda P'_C(x; v), \text{ for any } \lambda > 0.$$

4. *Let $y \in C$. Suppose $(P_C^{-1}(y))^o \neq \emptyset$. Then, $P_C$ is directionally differentiable on $(P_C^{-1}(y))^o$ such that, for any $x \in (P_C^{-1}(y))^o$, we have*

$$P'_C(x)(v) = \theta, \text{ for every } v \in Z \text{ with } v \neq \theta.$$

In particular, if $C$ is a closed ball, or closed and convex cone (including proper subspace) in $Z$, then the exact analytic representations of $P'_C$ are provided in [17]. In [18], the authors studied properties of $P'_C$ in Hilbert spaces and Hilbertian Bochner spaces, which are considered as special cases of uniformly convex and uniformly smooth Banach spaces.

Let $S(Z)$ denote the unit sphere in $Z$. We define the Fréchet differentiability of the metric projection.

**Definition 5.1**. Let $x \in Z$. Suppose that $P_C$ is (Gâteaux) directionally differentiable at point $x$. If the limit (5.1) is attained uniformly for $v \in S(Z)$, then $P_C$ is said to be Fréchet differentiable at point $x$. Let $D$ be a subset of $Z$. If $P_C$ is Fréchet differentiable at every point $x$ of $D$, then $P_C$ is said to be Fréchet differentiable on $D$.

## 5.2. The directional differentiability of the metric projection onto closed subspaces in uniformly convex and uniformly smooth Bochner spaces

The first theorem in this section proves the Fréchet differentiability of the metric projection operator onto closed subspaces of $L_p(S; X)$.

**Theorem 5.2.** *$P_{L_p(A;X)}$ is Fréchet differentiable on $L_p(S; X)$ such that, for any $f \in L_p(S; X)$, we have,*

$$P'_{L_p(A;X)}(f)(h) = h_A, \quad \text{for any } h \in L_p(S;X) \text{ with } h \neq \theta.$$

*Proof.* For any $f, h \in L_p(S;X)$ with $h \neq \theta$, by Lemma 4.1 and Lemma 3.3, we have

$$
\begin{aligned}
&P'_{L_p(A;X)}(f)(h) \\
&= \lim_{t \downarrow 0} \frac{P_{L_p(A;X)}(f+th) - P_{L_p(A;X)}(f)}{t} \\
&= \lim_{t \downarrow 0} \frac{(f+th)_A - f_A}{t} \\
&= \lim_{t \downarrow 0} \frac{f_A + t h_A - f_A}{t} \\
&= h_A.
\end{aligned}
$$

It is clear that the above limit uniformly converge on $S(L_p(S;X))$. □

**Remarks 5.3.** From Theorem 5.2, we see that $P'_{L_p(A;X)}$ has the following properties.

(a) The directional derivatives $P'_{L_p(A;X)}(f)(h)$ at point $f$ only depends on the direction $h$;

(b) For any $f, h \in L_p(S;X)$ with $h \neq \theta$, $P'_{L_p(A;X)}(f)(h) = h_A \in L_p(A;X)$;

(c) If $\mu(S \backslash A) = 0$, then $L_p(A;X)$ coincides with $L_p(S;X)$ and $P_{L_p(A;X)}$ coincides with $P_{L_p(S;X)}$, which is the identity mapping in $L_p(S;X)$. In this case, Theorem 5.2 is trivial.

### 5.3. The directional differentiability of the metric projection onto closed balls in subspaces of uniformly convex and uniformly smooth Bochner spaces

In next theorem, we prove the directional differentiability of the metric projection operator onto closed balls in $L_p(A;X)$.

**Theorem 5.4**. *For $r > 0$ and for $g, h \in L_p(S;X)$ with $h \neq \theta$, we have,*

(a) *If $g \in B_A^o(r)$, then*

$$P'_{B_A(r)}(g)(h) = \begin{cases} h, & \text{for } h \in L_p(A;X) \\ h_A, & \text{for } h \notin L_p(A;X) \end{cases};$$

(b) *If $g \in L_p(A;X) \backslash B_A(r)$, then*

$$P'_{B_A(r)}(g)(h) = \begin{cases} \dfrac{r}{\|g\|_{L_p(A;X)}} \left( h - \dfrac{\Psi_p(g,h)}{\|g\|_{L_p(A;X)}} g \right), & \text{for } h \in L_p(A;X) \\ \dfrac{r}{\|g\|_{L_p(A;X)}} \left( h_A - \dfrac{\Psi_p(g,h_A)}{\|g\|_{L_p(A;X)}} g \right), & \text{for } h \notin L_p(A;X) \end{cases}.$$

*In particular,*

$$P'_{B_A(r)}(g)(g) = \theta, \text{ for every } g \in L_p(A;X) \backslash B_A(r);$$

(c) *If $g \in C_A^o(r) \backslash L_p(A;X)$, then*

$$P'_{B_A(r)}(g)(h) = h_A;$$

(d) *If $g \in L_p(S;X) \backslash (C_A(r) \cup L_p(A;X))$, then*

$$P'_{B_A(r)}(g)(h) = \frac{r}{\|g_A\|_{L_p(A;X)}} \left( h_A - \frac{\Psi_p(g_A, h_A)}{\|g_A\|_{L_p(A;X)}} g_A \right).$$

*In particular,*

$$P'_{B_A(r)}(g)(g) = \theta, \text{ for every } g \in L_p(S;X) \backslash (C_A(r) \cup L_p(A;X)).$$

*Proof.* We prove (a). Suppose $g \in B_A^o(r)$. If $h \in L_p(A;X)$, then, there is $\delta > 0$, such that $g + th \in B_A(r)$, for all $t \in (0, \delta)$. By part (a) in Theorem 4.4, we have

$$\begin{aligned}
&P'_{B_A(r)}(g)(h) \\
&= \lim_{t \downarrow 0} \frac{P_{B_A(r)}(g+th) - P_{B_A(r)}(g)}{t} \\
&= \lim_{t \downarrow 0, t < \delta} \frac{g + th - g}{t} \\
&= \lim_{t \downarrow 0, t < \delta} \frac{th}{t} \\
&= h.
\end{aligned}$$

If $h \in L_p(S;X) \backslash L_p(A;X)$, then, there is $\lambda > 0$, such that $g + th \in C_A(r) \backslash L_p(A;X)$, for all $t \in (0, \lambda)$. By parts (a, c) in Theorem 4.4 and by Lamma 3.3, we have

$$\begin{aligned}
&P'_{B_A(r)}(g)(h) \\
&= \lim_{t \downarrow 0} \frac{P_{B_A(r)}(g+th) - P_{B_A(r)}(g)}{t} \\
&= \lim_{t \downarrow 0, t < \lambda} \frac{(g+th)_A - g}{t} \\
&= \lim_{t \downarrow 0, t < \lambda} \frac{g + th_A - g}{t} \\
&= \lim_{t \downarrow 0, t < \lambda} \frac{th_A}{t} \\
&= h_A.
\end{aligned}$$

Proof of (b). Let $g \in L_p(A;X) \backslash B_A(r)$. It implies $\|g\|_{L_p(A;X)} > r$. If $h \in L_p(A;X)$, then, there is $\delta > 0$, such that $g + th \in L_p(A;X) \backslash B_A(r)$, for all $t \in (0, \delta)$. By part (b) in Theorem 4.4, we have

$$\begin{aligned}
&P'_{B_A(r)}(g)(h) \\
&= \lim_{t \downarrow 0} \frac{P_{B_A(r)}(g+th) - P_{B_A(r)}(g)}{t} \\
&= \lim_{t \downarrow 0, t < \delta} \frac{\frac{r}{\|(g+th)\|_{L_p(A;X)}}(g+th) - \frac{r}{\|g\|_{L_p(A;X)}}g}{t} \\
&= \lim_{t \downarrow 0, t < \delta} \frac{\frac{rth}{\|(g+th)\|_{L_p(A;X)}} + \left( \frac{r}{\|(g+th)\|_{L_p(A;X)}} - \frac{r}{\|g\|_{L_p(A;X)}} \right)g}{t} \\
&= \lim_{t \downarrow 0, t < \delta} \frac{rh}{\|(g+th)\|_{L_p(A;X)}} + \lim_{t \downarrow 0, t < \delta} \frac{\left( \frac{r}{\|(g+th)\|_{L_p(A;X)}} - \frac{r}{\|g\|_{L_p(A;X)}} \right)g}{t} \\
&= \frac{r}{\|g\|_{L_p(A;X)}}h - \lim_{t \downarrow 0, t < \delta} \frac{\frac{r}{\|(g+th)\|_{L_p(A;X)}\|g\|_{L_p(A;X)}}\left( \|(g+th)\|_{L_p(A;X)} - \|g\|_{L_p(A;X)} \right)g}{t} \\
&= \frac{r}{\|g\|_{L_p(A;X)}}h - \frac{r}{\|g\|_{L_p(A;X)}^2}\Psi_p(g,h)g
\end{aligned}$$

$$= \frac{r}{\|g\|_{L_p(A;X)}} \left( h - \frac{\Psi_p(g,h)}{\|g\|_{L_p(A;X)}} g \right).$$

If $h \in L_p(S;X) \backslash L_p(A;X)$, by $g \in L_p(A;X) \backslash B_A(r)$ with $\|g\|_{L_p(A;X)} > r$, there is $\lambda > 0$, such that $g + th \in L_p(S;X) \backslash (C_A(r) \cup L_p(A;X))$, for all $t \in (0, \lambda)$. By parts (b), (d) in Theorem 4.4 and by Lamma 3.3, we have

$$
\begin{aligned}
& P'_{B_A(r)}(g)(h) \\
&= \lim_{t \downarrow 0} \frac{P_{B_A(r)}(g+th) - P_{B_A(r)}(g)}{t} \\
&= \lim_{t \downarrow 0, t < \lambda} \frac{\frac{r}{\|(g+th)_A\|_{L_p(A;X)}}(g+th)_A - \frac{r}{\|g\|_{L_p(A;X)}}g}{t} \\
&= \lim_{t \downarrow 0, t < \lambda} \frac{\frac{rth_A}{\|(g+th)_A\|_{L_p(A;X)}} + \left( \frac{r}{\|(g+th)_A\|_{L_p(A;X)}} - \frac{r}{\|g\|_{L_p(A;X)}} \right) g}{t} \\
&= \lim_{t \downarrow 0, t < \lambda} \frac{rh_A}{\|(g+th)_A\|_{L_p(A;X)}} - \lim_{t \downarrow 0, t < \lambda} \frac{\frac{r}{\|(g+th)_A\|_{L_p(A;X)}\|g\|_{L_p(A;X)}} \left( \|(g+th)_A\|_{L_p(A;X)} - \|g\|_{L_p(A;X)} \right) g}{t} \\
&= \lim_{t \downarrow 0, t < \lambda} \frac{rh_A}{\|(g+th)_A\|_{L_p(A;X)}} - \lim_{t \downarrow 0, t < \lambda} \frac{\frac{r}{\|(g+th)_A\|_{L_p(A;X)}\|g\|_{L_p(A;X)}} \left( \|g+th_A\|_{L_p(A;X)} - \|g\|_{L_p(A;X)} \right) g}{t} \\
&= \frac{r}{\|g\|_{L_p(A;X)}} h - \frac{r}{\|g\|_{L_p(A;X)}^2} \Psi_p(g, h_A) g \\
&= \frac{r}{\|g\|_{L_p(A;X)}} \left( h_A - \frac{\Psi_p(g,h_A)}{\|g\|_{L_p(A;X)}} g \right).
\end{aligned}
$$

Then, part (b) is completed proved by using Lemma 2.2.

Proof of (c). Let $g \in C_A^o(r) \backslash L_p(A;X)$ with $\|g_A\|_{L_p(A;X)} < r$. For $h \in L_p(S;X)$ with $h \neq \theta$, there is $\delta > 0$, such that $g + th \in C_A^o(r) \backslash L_p(A;X)$, for all $t \in (0, \delta)$. By part (c) in Theorem 4.4 and by Lamma 3.3, we have

$$
\begin{aligned}
& P'_{B_A(r)}(g)(h) \\
&= \lim_{t \downarrow 0} \frac{P_{B_A(r)}(g+th) - P_{B_A(r)}(g)}{t} \\
&= \lim_{t \downarrow 0, t < \delta} \frac{(g+th)_A - g_A}{t} \\
&= \lim_{t \downarrow 0, t < \delta} \frac{th_A}{t} \\
&= h_A.
\end{aligned}
$$

Proof of (d). Let $g \in L_p(S;X) \backslash (C_A(r) \cup L_p(A;X))$ with $\|g_A\|_{L_p(A;X)} > r$. For $h \in L_p(S;X)$ with $h \neq \theta$, there is $\delta > 0$, such that $g + th \notin L_p(A;X))$ and $\|(g+th)_A\|_{L_p(A;X)} > r$, that is, $g + th \in L_p(S; X) \backslash (C_A(r) \cup L_p(A;X))$, for all $t \in (0, \delta)$. By part (d) in Theorem 4.4 and by Lamma 3.3, we have

$$
\begin{aligned}
& P'_{B_A(r)}(g)(h) \\
&= \lim_{t \downarrow 0} \frac{P_{B_A(r)}(g+th) - P_{B_A(r)}(g)}{t} \\
&= \lim_{t \downarrow 0, t < \delta} \frac{\frac{r}{\|(g+th)_A\|_{L_p(A;X)}}(g+th)_A - \frac{r}{\|g_A\|_{L_p(A;X)}}g_A}{t}
\end{aligned}
$$

$$= \lim_{t \downarrow 0, t < \delta} \frac{\frac{r t h_A}{\|(g+th)A\|_{L_p(A;X)}} + \left(\frac{r}{\|(g+th)A\|_{L_p(A;X)}} - \frac{r}{\|g_A\|_{L_p(A;X)}}\right) g_A}{t}$$

$$= \lim_{t \downarrow 0, t < \delta} \frac{r h_A}{\|(g+th)A\|_{L_p(A;X)}} + \lim_{t \downarrow 0, t < \delta} \frac{\left(\frac{r}{\|(g+th)A\|_{L_p(A;X)}} - \frac{r}{\|g_A\|_{L_p(A;X)}}\right) g_A}{t}$$

$$= \frac{r}{\|g_A\|_{L_p(A;X)}} h_A - \lim_{t \downarrow 0, t < \delta} \frac{\frac{r}{\|(g+th)A\|_{L_p(A;X)}\|(g+th)A\|_{L_p(A;X)}} \left(\|(g_A+th_A\|_{L_p(A;X)} - \|g_A\|_{L_p(A;X)}\right) g_A}{t}$$

$$= \frac{r}{\|g_A\|_{L_p(A;X)}} h_A - \frac{r}{\|g_A\|^2_{L_p(A;X)}} \Psi_p(g_A, h_A) g_A$$

$$= \frac{r}{\|g_A\|_{L_p(A;X)}} \left(h_A - \frac{\Psi_p(g_A, h_A)}{\|g_A\|_{L_p(A;X)}} g_A\right).$$

Then, part (d) is completely proved by using Lemma 2.2. □

Notice again, if $\mu(S\backslash A) = 0$, then $L_p(A; X)$ coincides with $L_p(S; X)$. In this case, we have

$$B_A(r) = B(r) \quad \text{and} \quad C_A(r) = B(r).$$

Furthermore, we have

$$h_A = h, \text{ for any } h \in L_p(S; X).$$

Then, by parts (a) and (b) in Theorem 5.4, we have

**Corollary 5.5**. *For $r > 0$ and for $g, h \in L_p(S; X)$ with $h \neq \theta$, we have,*

(a) *If $g \in B^o(r)$, then*

$$P'_{B(r)}(g)(h) = h, \text{for any } h \in L_p(S; X) \text{ with } h \neq \theta;$$

(b) *If $g \in L_p(S; X)\backslash B(r)$, then*

$$P'_{B(r)}(g)(h) = \frac{r}{\|g\|_{L_p(S;X)}} \left(h - \frac{\Psi_p(g,h)}{\|g\|_{L_p(S;X)}} g\right), \text{for any } h \in L_p(S; X) \text{ with } h \neq \theta.$$

*In particular,*

$$P'_{B(r)}(g)(g) = \theta, \text{for every } g \in L_p(S; X)\backslash B(r).$$

## 5.4. The directional differentiability of the metric projection onto closed and convex cylinders in uniformly convex and uniformly smooth Bochner spaces

In the following theorem, we prove the Fréchet differentiability of the metric projection operator onto closed and convex cylinders in $L_p(S; X)$.

**Theorem 5.6**. *For $r > 0$, we have,*

(a) *$P_{C_A(r)}$ is Fréchet differentiable on $C_A^o(r)$ satisfying that, for any $g \in C_A^o(r)$,*

$$P'_{C_A(r)}(g)(h) = h, \text{ for any } h \in L_p(S; X)\backslash\{\theta\};$$

(b) *$P_{C_A(r)}$ is Fréchet differentiable on $L_p(S; X)\backslash C_A(r)$ such that, for $g \in L_p(S; X)\backslash C_A(r)$,*

$$P'_{C_A(r)}(g)(h) = \frac{r}{\|g_A\|_{L_p(A;X)}}\left(h_A + h_{S\setminus A} - \frac{\Psi_p(g_A, h_A)}{\|g_A\|_{L_p(A;X)}}g_A\right), \textit{for any } h \in L_p(S;X)\setminus\{\theta\}.$$

*In particular,*

$$P'_{C_A(r)}(g)(g) = g_{S\setminus A}, \textit{for every } g \in L_p(S;X)\setminus C_A(r).$$

*Proof.* We prove (a). Let $g \in C_A^o(r)$. Since $\|g_A\|_{L_p(A;X)} < r$, then, for any $h \in L_p(S;X)$ with $h \neq \theta$, there is $\delta > 0$, such that $\|(g+th)_A\|_{L_p(A;X)} < r$, that is, $g + th \in C_A^o(r)$, for all $t \in (0, \delta)$. By part (a) in Theorem 4.7, we have

$$\begin{aligned}
&P'_{C_A(r)}(g)(h)\\
&= \lim_{t\downarrow 0}\frac{P_{C_A(r)}(g+th) - P_{C_A(r)}(g)}{t}\\
&= \lim_{t\downarrow 0, t<\delta}\frac{g+th-g}{t}\\
&= h.
\end{aligned}$$

Notice that, for a given $g \in C_A^o(r)$, we can choose a positive $\delta$ with $\delta < \frac{r - \|g_A\|_{L_p(A;X)}}{2}$ to satisfy $g + th \in C_A^o(r)$, for all $t \in (0, \delta)$, which is independent from $h$ with $\|h\|_{L_p(S;X)} = 1$. This implies that the above limit is uniformly convergent with respect to $h$ satisfying $\|h\|_{L_p(S;X)} = 1$. This implies the Fréchet differentiability of the metric projection operator $P_{C_A(r)}$ onto closed and convex cylinder $C_A(r)$ in $L_p(S;X)$ at the point $g \in C_A^o(r)$.

Proof of (b). Let $g \in L_p(S;X)\setminus C_A(r)$. Since $\|g_A\|_{L_p(A;X)} > r$, then, we can choose a number $\delta$ with $0 < \delta < \frac{\|g_A\|_{L_p(A;X)} - r}{2}$, such that, for any $t \in (0, \delta)$ and any $h \in L_p(S;X)$ with $\|h\|_{L_p(S;X)} = 1$, we have

$$\|(g+th)_A\|_{L_p(A;X)} \geq \|g_A\|_{L_p(A;X)} - \frac{\|g_A\|_{L_p(A;X)} - r}{2} > r.$$

This implies

$$g + th \in L_p(S;X)\setminus C_A(r), \text{ for any } t \in (0, \delta) \text{ and any } h \in L_p(S;X) \text{ with } \|h\|_{L_p(S;X)} = 1.$$

By part (b) in Theorem 4.7 and Lemma 3.3, we have

$$\begin{aligned}
&P'_{C_A(r)}(g)(h)\\
&= \lim_{t\downarrow 0}\frac{P_{C_A(r)}(g+th) - P_{C_A(r)}(g)}{t}\\
&= \lim_{t\downarrow 0, t<\delta}\frac{\left(\frac{r}{\|(g+th)_A\|_{L_p(A;X)}}(g+th)_A + (g+th)_{S\setminus A}\right) - \left(\frac{r}{\|g_A\|_{L_p(A;X)}}g_A + g_{S\setminus A}\right)}{t}\\
&= \lim_{t\downarrow 0, t<\delta}\frac{\frac{r}{\|(g+th)_A\|_{L_p(A;X)}}th_A + th_{S\setminus A} + \left(\frac{r}{\|(g+th)_A\|_{L_p(A;X)}} - \frac{r}{\|g_A\|_{L_p(A;X)}}\right)g_A}{t}\\
&= \frac{r}{\|g_A\|_{L_p(A;X)}}h_A + h_{S\setminus A} + \lim_{t\downarrow 0, t<\delta}\frac{\left(\frac{r}{\|(g+th)_A\|_{L_p(A;X)}} - \frac{r}{\|g_A\|_{L_p(A;X)}}\right)g_A}{t}\\
&= \frac{r}{\|g_A\|_{L_p(A;X)}}h_A + h_{S\setminus A} - \lim_{t\downarrow 0, t<\delta}\frac{\frac{r}{\|(g+th)_A\|_{L_p(A;X)}\|(g+th)_A\|_{L_p(A;X)}}\left(\|g_A+th_A\|_{L_p(A;X)} - \|g_A\|_{L_p(A;X)}\right)g_A}{t}
\end{aligned}$$

$$= \frac{r}{\|g_A\|_{L_p(A;\,X)}} h_A + h_{S\setminus A} - \frac{r}{\|g_A\|_{L_p(A;\,X)}^2} \Psi_p(g_A, h_A) g_A$$

$$= \frac{r}{\|g_A\|_{L_p(A;\,X)}} \left( h_A + h_{S\setminus A} - \frac{\Psi_p(g_A, h_A)}{\|g_A\|_{L_p(A;\,X)}} g_A \right).$$

It is clear that the above limit is uniformly convergent with respect to $h \in L_p(S; X)$ satisfying $\|h\|_{L_p(S;\,X)}$ = 1. Then, part (b) is completely proved by using Lemma 2.2. $\qquad\square$

**Remarks on Corollary 5.5 and Theorem 5.6.** As we noticed before Corollary 4.5, if $\mu(S\setminus A) = 0$, then $L_p(A; X)$ coincides with $L_p(S; X)$. In this case, we have $C_A(r) = B(r)$ and

$$h_{S\setminus A} = \theta \text{ and } h_A = h, \text{ for } h \in L_p(S; X).$$

Then, we see that Corollary 5.5 follows from Theorem 5.6 immediately.

## 6. Applications to Hilbert spaces with orthonormal bases

### 6.1. The metric projection onto closed balls in Hilbert spaces

In this section, we study the directional differentiability of the metric projection in Hilbert spaces, which are considered as special cases of uniformly convex and uniformly smooth Banach spaces. Or, they can be considered as special cases of Hilbertian Bochner spaces.

Through this section, as usual, let $\mathbb{N}$ denote the set of all positive integers. Let $(H, \|\cdot\|)$ be a real Hilbert space with inner product $\langle \cdot, \cdot \rangle$. Suppose that $H$ has an orthonormal basis $\{e_n\}_{n\in N}$, in which $N$ is a subset of $\mathbb{N}$ with cardinality greater than 1 and $\{e_n\}_{n\in N}$ satisfies that, for $m, n \in N$,

$$\langle e_m, e_n \rangle = \begin{cases} 1, & \text{if } m = n, \\ 0, & \text{if } m \neq n. \end{cases}$$

Every $x \in H$ has the following analytic representation

$$x = \sum_{n\in N} \langle x, e_n \rangle e_n \quad \text{such that} \quad \|x\|^2 = \sum_{n\in N} \langle x, e_n \rangle^2,$$

and

$$\langle x, y \rangle = \sum_{n\in N} \langle x, e_n \rangle \langle y, e_n \rangle, \text{ for any } x, y \in H.$$

Let $M$ be a nonempty subset of $N$. Let $H_M$ be the closed subspace of $H$ generated by $\{e_m\}_{m\in M}$. If $N\setminus M \neq \emptyset$, then $H_M$ is a closed proper subspace of $H$. $H_{N\setminus M}$ is similarly defined to be the closed subspace of $H$ generated by $\{e_n\}_{n\in N\setminus M}$ that satisfies

$$(H_M)^{\perp} = H_{N\setminus M}. \tag{6.1}$$

That is, $H_M$ and $H_{N\setminus M}$ are orthogonal spaces of each other. For any $x \in H$, we define $x_M \in H_M$ by

$$x_M = \sum_{m\in M} \langle x, e_m \rangle e_m.$$

For any $r > 0$, the closed, open balls and the sphere in $H_M$ with radius $r$ and with center at the origin are respectively denoted by

$$B_M(r) = \left\{ x \in H_M : \left( \sum_{m\in M} \langle x, e_m \rangle^2 \right)^{\frac{1}{2}} \leq r \right\},$$

$$B_M^o(r) = \left\{ x \in H_M : \ (\textstyle\sum_{m \in M} \langle x, e_m \rangle^2)^{\frac{1}{2}} < r \right\},$$
$$S_M(r) = \left\{ x \in H_M : \ (\textstyle\sum_{m \in M} \langle x, e_m \rangle^2)^{\frac{1}{2}} = r \right\}.$$

$B_M(r)$ is a nonempty closed, bounded and convex subset in $H$. If $N \backslash M \neq \emptyset$, then $B_M(r)$ is not a closed ball in $H$. We define the following subsets in $H$:

$$C_M(r) = \{ x \in H : x_M \in B_M(r) \},$$

$$C_M^o(r) = \{ x \in H : x_M \in B_M^o(r) \}.$$

$C_M(r)$ and $C_M^o(r)$ are the convex cylinders in $H$ with bases $B_M(r)$ and $B_M^o(r)$ in $H_M$, respectively. If $N \backslash M \neq \emptyset$, then $C_M(r)$ is a closed, unbounded and convex subset in $H$ and $C_A^o(r)$ is an open, unbounded and convex subset in $H$. In particular, if $N \backslash M = \emptyset$, then $H_M$ coincides with $H$ and, in this case, $B_M(r)$, $B_M^o(r)$, $S_M(r)$, $C_M(r)$ and $C_M^o(r)$ are respectively denoted by $B(r)$, $B^o(r)$, $S(r)$, $C(r)$ and $C^o(r)$, which satisfy

$$B(r) = C(r) \quad \text{and} \quad B^o(r) = C^o(r).$$

Next, we calculate the values of the metric projection from $H$ onto $B_M(r)$.

**Theorem 6.1.** *For any $r > 0$ and $x \in H$, we have*

(a) $P_{B_M(r)}(x) = x$, *for $x \in B_M(r)$;*
(b) $P_{B_M(r)}(x) = \frac{r}{\|x\|} x$, *for $x \in H_M \backslash B_M(r)$;*
(c) $P_{B_M(r)}(x) = x_M$, *for $x \in C_M(r) \backslash H_M$;*
(d) $P_{B_M(r)}(x) = \frac{r}{\|x_M\|} x_M$, *for $x \in H \backslash (C_M(r) \cup H_M)$.*

*Proof.* The proof of Theorem 6.1 is similar to the proof of Theorem 4.4. Part (a) is clear. We prove part (b). For any given $x \in H_M \backslash B_M(r)$, $x$ must satisfies $\|x\| > r$. Then, for any $y \in B_M(r)$,

$$\begin{aligned}
& \left\langle x - \frac{r}{\|x\|} x, \ \frac{r}{\|x\|} x - y \right\rangle \\
&= \left( 1 - \frac{r}{\|x\|} \right) \left\langle x, \ \frac{r}{\|x\|} x - y \right\rangle \\
&= \left( 1 - \frac{r}{\|x\|} \right) \left( r \|x\| - \langle x, y \rangle \right) \\
&\geq (\|x\| - r)\,(r - \|y\|) \\
&\geq 0, \text{ for all } y \in B_M(r).
\end{aligned}$$

By $x \in H_M$ and $\|x\| > r$, we have $\frac{r}{\|x\|} x \in S_M(r)$. By the basic variational principle of $P_{B_M(r)}$ in Hilbert spaces, the above inequality implies

$$P_{B_M(r)}(x) = \frac{r}{\|x\|} x, \quad \text{for any } x \in H_M \backslash B_M(r).$$

Proof of (c). For $x \in C_M(r) \backslash H_M$, we have $x_M \in B_M(r)$ and $x \notin H_M$. It follows that

$$\|x_M\| \leq r \quad \text{and} \quad 0 < \|x_M\| < \|x\|. \tag{6.2}$$

Then, for any $y \in B_M(r) \subseteq H_M$, we have $x_M - y \in H_M$. By (6.1), we obtain

$$\langle x - x_M, \ x_M - y \rangle$$
$$= \langle x_{N\backslash M}, x_M - y \rangle$$
$$= 0, \text{ for all } y \in B_M(r).$$

By (6.2) and the basic variational principle of $P_{B_M(r)}$, this implies $P_{B_M(r)}(x) = x_M$.

The proof of part (d) is similar to the proof of part (b). For $x \in H\backslash(C_M(r) \cup H_M)$, we have $x_M \notin B_M(r)$ and $x \notin H_M$, it follows that

$$r < \|x_M\| < \|x\|. \tag{6.3}$$

Then, for any $y \in B_M(r) \subseteq H_M$ and by (6.1) and (6.3), we have

$$\langle x - \frac{r}{\|x_M\|} x_M, \ \frac{r}{\|x_M\|} x_M - y \rangle$$
$$= \langle x_{N\backslash M} + \left(1 - \frac{r}{\|x_M\|}\right) x_M, \ \frac{r}{\|x_M\|} x_M - y \rangle$$
$$= \left(1 - \frac{r}{\|x_M\|}\right) \langle x_M, \ \frac{r}{\|x_M\|} x_M - y \rangle$$
$$= \left(1 - \frac{r}{\|x_M\|}\right) \left(r\|x_M\| - \langle x_M, y \rangle\right)$$
$$\geq \left(\|x_M\| - r\right) \left(r - \|y\|\right)$$
$$\geq 0, \text{ for all } y \in B_M(r).$$

By $\|x_M\| > r$, we have $\frac{r}{\|x_M\|} x_M \in S_M(r)$. By the basic variational principle of $P_{B_M(r)}$, this implies

$$P_{B_M(r)}(x) = \frac{r}{\|x_M\|} x_M, \quad \text{for } x \in H\backslash(C_M(r) \cup H_M). \qquad \square$$

Notice that if $N\backslash M = \emptyset$, then $H_M$ coincides with $H$. In this case, we have

$$B_M(r) = B(r) \quad \text{and} \quad C_M(r) = B(r).$$

Then, by Theorem 6.1, we have

**Corollary 6.2.** *For any $r > 0$ and $x \in H$, we have*

(a) $P_{B(r)}(x) = x, \ \text{for } x \in B(r);$

(b) $P_{B(r)}(x) = \frac{r}{\|x\|} x, \ \text{for } x \in H\backslash B(r).$

By Theorem 6.1, we can study the directional differentiability of $P_{B_M(r)}$.

**Theorem 6.3.** *For any $r > 0$ and $x, h \in H$ with $h \neq \theta$, we have,*

(a) *If $x \in B_M^o(r)$, then*

$$P'_{B_M(r)}(x)(h) = \begin{cases} h, & \text{for } h \in H_M \\ h_M, & \text{for } h \notin H_M \end{cases};$$

(b) *If $x \in H_M\backslash B_M(r)$, then*

$$P'_{B_M(r)}(x)(h) = \begin{cases} \frac{r}{\|x\|}\left(h - \frac{\langle x, h \rangle}{\|x\|^2} x\right), & \text{for } h \in H_M \\ \frac{r}{\|x\|}\left(h_M - \frac{\langle x, h_M \rangle}{\|x\|^2} x\right), & \text{for } h \notin H_M \end{cases}.$$

*In particular,*

$$P'_{B_M(r)}(x)(x) = \theta, \text{ for any } x \in H_M \backslash B_M(r);$$

(c) *If $x \in C_M^o(r) \backslash H_M$, then*

$$P'_{B_M(r)}(x)(h) = h_M;$$

(d) *If $x \in H \backslash (C_M(r) \cup H_M)$, then*

$$P'_{B_M(r)}(x)(h) = \frac{r}{\|x_M\|} \left( h_M - \frac{\langle x_M, h_M \rangle}{\|x_M\|^2} x_M \right).$$

*In particular,*

$$P'_{B_M(r)}(x)(x) = \theta, \text{ for any } x \in H \backslash (C_M(r) \cup H_M).$$

*Proof.* The proof of this theorem is similar to the proof of Theorem 5.4. Proof of (a). Let $x \in B_M^o(r)$. If $h \in H_M$, then, there is $\delta > 0$, such that $x + th \in B_M(r)$, for all $t \in (0, \delta)$. By part (a) in Theorem 6.1, we have

$$\begin{aligned}
&P'_{B_M(r)}(x)(h) \\
&= \lim_{t \downarrow 0} \frac{P_{B_M(r)}(x+th) - P_{B_M(r)}(x)}{t} \\
&= \lim_{t \downarrow 0, t < \delta} \frac{x+th-x}{t} \\
&= \lim_{t \downarrow 0, t < \delta} \frac{th}{t} \\
&= h.
\end{aligned}$$

If $h \in H \backslash H_M$, then, there is $\lambda > 0$, such that $x + th \in C_M(r) \backslash H_M$, for all $t \in (0, \lambda)$. By parts (a, c) in Theorem 6.1 and by Lamma 3.3, we have

$$\begin{aligned}
&P'_{B_M(r)}(x)(h) \\
&= \lim_{t \downarrow 0} \frac{P_{B_M(r)}(x+th) - P_{B_M(r)}(x)}{t} \\
&= \lim_{t \downarrow 0, t < \lambda} \frac{(x+th)_M - x}{t} \\
&= \lim_{t \downarrow 0, t < \lambda} \frac{x+th_M - x}{t} \\
&= \lim_{t \downarrow 0, t < \lambda} \frac{th_M}{t} \\
&= h_M.
\end{aligned}$$

Proof of (b). Let $x \in H_M \backslash B_M(r)$. It implies $\|x\| > r$. If $h \in H_M$ with $h \neq \theta$, then, there is $\delta > 0$, such that $x + th \in H_M \backslash B_M(r)$, for all $t \in (0, \delta)$. Noticing $x \neq \theta$, by part (b) in Theorem 6.1 and (2.4) for the solutions of $\Psi_H$, we have

$$\begin{aligned}
&P'_{B_M(r)}(x)(h) \\
&= \lim_{t \downarrow 0} \frac{P_{B_M(r)}(x+th) - P_{B_M(r)}(x)}{t} \\
&= \lim_{t \downarrow 0, t < \delta} \frac{\frac{r}{\|x+th\|}(x+th) - \frac{r}{\|x\|}x}{t} \\
&= \lim_{t \downarrow 0, t < \delta} \frac{\frac{rth}{\|x+th\|} + \left(\frac{r}{\|x+th\|} - \frac{r}{\|x\|}\right)x}{t} \\
&= \lim_{t \downarrow 0, t < \delta} \frac{rh}{\|x+th\|} + \lim_{t \downarrow 0, t < \delta} \frac{\left(\frac{r}{\|x+th\|} - \frac{r}{\|x\|}\right)x}{t} \\
&= \frac{r}{\|x\|}h - \lim_{t \downarrow 0, t < \delta} \frac{\frac{r}{\|x+th\|\|x\|}(\|x+th\| - \|x\|)x}{t}
\end{aligned}$$

$$= \frac{r}{\|x\|} h - \frac{r}{\|x\|^2} \frac{1}{\|x\|} \langle x, h \rangle x$$
$$= \frac{r}{\|x\|} \left( h - \frac{\langle x,h \rangle}{\|x\|^2} x \right).$$

If $h \in H \backslash H_M$, by $x \in H_M \backslash B_M(r)$ with $\|x\| > r$, there is $\lambda > 0$, such that $x + th \in H \backslash (C_M(r) \cup H_M)$, for all $t \in (0, \lambda)$. Noticing $x \neq \theta$, by parts (b, d) in Theorem 6.1, Lamma 3.3 and (2.4), we have

$$P'_{B_M(r)}(x)(h)$$
$$= \lim_{t \downarrow 0} \frac{P_{B_M(r)}(x+th) - P_{B_M(r)}(x)}{t}$$
$$= \lim_{t \downarrow 0, t < \lambda} \frac{\frac{r}{\|(x+th)_M\|}(x+th)_M - \frac{r}{\|x\|} x}{t}$$
$$= \lim_{t \downarrow 0, t < \lambda} \frac{\frac{r(th)_M}{\|(x+th)_M\|} + \left( \frac{r}{\|(x+th)_M\|} - \frac{r}{\|x\|} \right) x}{t}$$
$$= \lim_{t \downarrow 0, t < \lambda} \frac{r h_M}{\|(x+th)_M\|} + \lim_{t \downarrow 0, t < \lambda} \frac{\left( \frac{r}{\|(x+th)_M\|} - \frac{r}{\|x\|} \right) x}{t}$$
$$= \frac{r}{\|x\|} h_M - \frac{r}{\|x\|^2} \frac{1}{\|x\|} \langle x, h_M \rangle x$$
$$= \frac{r}{\|x\|} \left( h_M - \frac{\langle x, h_M \rangle}{\|x\|^2} x \right).$$

Proof of (c). Let $x \in C_M^o(r) \backslash H_M$ with $\|x_M\| < r$. For $h \in H$ with $h \neq \theta$, there is $\delta > 0$, such that $x + th \in C_M^o(r) \backslash H_M$, for all $t \in (0, \delta)$. By part (c) in Theorem 6.1 and by Lamma 3.3, we have

$$P'_{B_M(r)}(x)(h)$$
$$= \lim_{t \downarrow 0} \frac{P_{B_M(r)}(x+th) - P_{B_M(r)}(x)}{t}$$
$$= \lim_{t \downarrow 0, t < \delta} \frac{(x+th)_M - x_M}{t}$$
$$= \lim_{t \downarrow 0, t < \delta} \frac{th_M}{t}$$
$$= h_M.$$

Proof of (d). Let $x \in H \backslash (C_M(r) \cup H_M)$ with $\|x\| > \|x_M\| > r$. For $h \in H$ with $h \neq \theta$, there is $\delta > 0$, such that $\|x + th\| > \|(x + th)_M\| > r$, that is, $x + th \in H \backslash (C_M(r) \cup H_M)$, for all $t \in (0, \delta)$. By part (d) in Theorem 6.1 and by Lamma 3.3, we have

$$P'_{B_M(r)}(x)(h)$$
$$= \lim_{t \downarrow 0} \frac{P_{B_M(r)}(x+th) - P_{B_M(r)}(x)}{t}$$
$$= \lim_{t \downarrow 0, t < \delta} \frac{\frac{r}{\|(x+th)_M\|}(x+th)_M - \frac{r}{\|x_M\|} x_M}{t}$$
$$= \lim_{t \downarrow 0, t < \delta} \frac{\frac{r th_M}{\|(x+th)_M\|} + \left( \frac{r}{\|(x+th)_M\|} - \frac{r}{\|x_M\|} \right) x_M}{t}$$
$$= \lim_{t \downarrow 0, t < \delta} \frac{r h_M}{\|(x+th)_M\|} + \lim_{t \downarrow 0, t < \delta} \frac{\left( \frac{r}{\|(x+th)_M\|} - \frac{r}{\|x_M\|} \right) x_M}{t}$$
$$= \frac{r}{\|x_M\|} h_M - \frac{r}{\|x_M\|^2} \frac{1}{\|x_M\|} \langle x_M, h_M \rangle x_M$$
$$= \frac{r}{\|x_M\|} \left( h_M - \frac{\langle x_M, h_M \rangle}{\|x_M\|^2} x_M \right). \qquad \square$$

If $N \backslash M = \emptyset$, then $H_M$ coincides with $H$. In this case, we have

$$B_M(r) = C_M(r) = B(r).$$

Furthermore, we have

$$h_M = h, \text{ for any } h \in H.$$

Then, by Theorem 6.3, we have

**Corollary 6.4**. *Let $r > 0$. For $x, h \in H$ with $h \neq \theta$, we have,*

(a) *If $x \in B^o(r)$, then*

$$P'_{B(r)}(x)(h) = h, \text{for any } h \in H \text{ with } h \neq \theta;$$

(b) *If $x \in H \backslash B(r)$, then*

$$P'_{B(r)}(x)(h) = \frac{r}{\|x\|}\left(h - \frac{\langle x, h \rangle}{\|x\|^2}x\right), \text{ for any } h \in H \text{ with } h \neq \theta.$$

*In particular,*

$$P'_{B(r)}(x)(x) = \theta, \text{for any } x \in H \backslash B(r).$$

### 6.2. The metric projection onto closed and convex cylinders in Hilbert spaces

Similar to section 4, in this subsection, we first calculate the values of the metric projection from $H$ onto a closed and convex cylinder $C_M(r)$ in $H$.

**Theorem 6.5**. *For any $r > 0$ and $x \in H$, we have*

(a) $P_{C_M(r)}(x) = x, \text{ for } x \in C_M(r);$
(b) $P_{C_M(r)}(x) = \frac{r}{\|x_M\|}x_M + x_{N \backslash M}, \text{ for } x \in H \backslash C_M(r).$

*Proof.* The proof of this theorem is similar to the proof of Theorem 4.7. Part (a) is clear. We prove (b). For any given $x \in H \backslash C_M(r)$, $x$ must satisfies $\|x\| \geq \|x_M\| > r$. Then, for any $y \in C_M(r)$ satisfying $\|y_M\| \leq r$, we have

$$\langle x - \left(\frac{r}{\|x_M\|}x_M + x_{N \backslash M}\right), \frac{r}{\|x_M\|}x_M + x_{N \backslash M} - y\rangle$$
$$= \left(1 - \frac{r}{\|x_M\|}\right)\langle x_M, \frac{r}{\|x_M\|}x_M + x_{N \backslash M} - (y_M + y_{N \backslash M})\rangle$$
$$= \left(1 - \frac{r}{\|x_M\|}\right)(r\|x_M\| - \langle x_M, y_M\rangle)$$
$$\geq (\|x_M\| - r)(r - \|y_M\|)$$
$$\geq 0, \text{ for all } y \in C_M(r).$$

By $x \in H$ and $\|x_M\| > r$, we have $\frac{r}{\|x_M\|}x_M \in S_M(r)$. It implies that $\frac{r}{\|x_M\|}x_M \in C_M(r)$. By the basic variational principle of $P_{C_M(r)}$, the above inequality implies

$$P_{C_M(r)}(x) = \frac{r}{\|x_M\|}x_M, \text{ for any } x \in H \backslash C_M(r). \qquad \square$$

Now, we prove the Fréchet differentiability and calculate the directional derivatives of the metric projection operator $P$ onto closed and convex cylinders in $H$ in the following theorem.

**Theorem 6.6**. *For any $r > 0$, we have,*

(a) $P_{C_M(r)}$ is Fréchet differentiable on $C_M^o(r)$ satisfying that, for any $x \in C_M^o(r)$,

$$P'_{C_M(r)}(x)(h) = h, \text{ for any } h \in H \backslash \{\theta\};$$

(b) $P_{C_M(r)}$ is Fréchet differentiable on $H \backslash C_M(r)$ such that, for $x \in H \backslash C_M(r)$,

$$P'_{C_M(r)}(x)(h) = \frac{r}{\|x_M\|}\left(h_M - \frac{\langle x_M, h_M \rangle}{\|x_M\|^2} x_M\right) + h_{N \backslash M}, \text{ for any } h \in H \backslash \{\theta\}.$$

*In particular,*

$$P'_{C_M(r)}(x)(x) = x_{N \backslash M}, \text{ for any } x \in H \backslash C_M(r).$$

*Proof.* Proof of (a). Let $x \in C_M^o(r)$. Since $\|x_M\| < r$, then, for any $h \in H$ with $h \neq \theta$, there is $\delta > 0$, such that $\|(x + th)_M\| < r$, that is, $x + th \in C_M^o(r)$, for all $t \in (0, \delta)$. In fact, we can choose a positive $\delta$ with $\delta < \frac{r - \|x_M\|}{2}$ to satisfy $x + th \in C_M^o(r)$, for all $t \in (0, \delta)$, which is independent from $h$ with $\|h\| = 1$. By part (a) in Theorem 6.5, we have

$$P'_{C_M(r)}(x)(h)$$
$$= \lim_{t \downarrow 0} \frac{P_{C_M(r)}(x+th) - P_{C_M(r)}(x)}{t}$$
$$= \lim_{t \downarrow 0, t < \delta} \frac{x + th - x}{t}$$
$$= \lim_{t \downarrow 0, t < \delta} \frac{th}{t}$$
$$= h.$$

Notice that, for this given $x \in C_M^o(r)$, a positive $\delta$ is chosen such that it is independent from $h$ with $\|h\| = 1$. This implies that the above limit is uniformly convergent with respect to $h$ satisfying $\|h\| = 1$. This implies the Fréchet differentiability of the metric projection operator $P_{C_M(r)}$ onto closed and convex cylinder $C_M(r)$ of $H$ at the point $x \in C_M^o(r)$.

Proof of (b). Let $x \in H \backslash C_M(r)$. Since $\|x_M\| > r$, then, we can choose a number $\lambda$ with $0 < \lambda < \frac{\|x_M\| - r}{2}$, such that, for any $t \in (0, \lambda)$ and any $h \in H$ with $\|h\| = 1$, we have

$$\|(x + th)_M\| \geq \|x_M\| - \frac{\|x_M\| - r}{2} > r.$$

This implies

$$x + th \in H \backslash C_M(r), \text{ for any } t \in (0, \lambda) \text{ and any } h \in H \text{ with } \|h\| = 1.$$

By part (b) in Theorem 6.5 and Lemma 3.3, we have

$$P'_{C_M(r)}(x)(h)$$
$$= \lim_{t \downarrow 0} \frac{P_{C_M(r)}(x+th) - P_{C_M(r)}(x)}{t}$$
$$= \lim_{t \downarrow 0, t < \lambda} \frac{\left(\frac{r}{\|(x+th)_M\|}(x+th)_M + (x+th)_{N \backslash M}\right) - \left(\frac{r}{\|x_M\|}x_M + x_{N \backslash M}\right)}{t}$$
$$= \lim_{t \downarrow 0, t < \lambda} \frac{\frac{r(th)_M}{\|(x+th)_M\|} + \left(\frac{r}{\|(x+th)_M\|} - \frac{r}{\|x_M\|}\right)x_M + (th)_{N \backslash M}}{t}$$
$$= \lim_{t \downarrow 0, t < \lambda} \frac{rh_M}{\|(x+th)_M\|} + h_{N \backslash M} + \lim_{t \downarrow 0, t < \lambda} \frac{\left(\frac{r}{\|(x+th)_M\|} - \frac{r}{\|x_M\|}\right)x_M}{t}$$

$$= \frac{r}{\|x_M\|}h_M + h_{N\setminus M} - \frac{r}{\|x_M\|^2}\frac{1}{\|x_M\|}\langle x_M, h_M\rangle x_M$$

$$= \frac{r}{\|x_M\|}\left(h_M - \frac{\langle x_M, h_M\rangle}{\|x_M\|^2}x_M\right) + h_{N\setminus M}. \qquad \square$$